\documentclass{amsart}

\makeatletter
\@namedef{subjclassname@2020}{\textup{2020} Mathematics Subject Classification}
\makeatother

\usepackage{graphicx}%
\usepackage{multirow}%
\usepackage{amsmath,amssymb,amsfonts}%
\usepackage{amsthm}%
\usepackage{mathrsfs}%
\usepackage[title]{appendix}%
\usepackage{xcolor}%
\usepackage{textcomp}%
\usepackage{manyfoot}%
\usepackage{booktabs}%
\usepackage{algorithm}%
\usepackage{algorithmicx}%
\usepackage{algpseudocode}%
\usepackage{listings}%
\usepackage{colonequals}

\usepackage{amstext,amscd,amsbsy,amsxtra,latexsym}
\usepackage{amssymb}

\usepackage{hyperref}
\hypersetup{
	colorlinks   = true, %Colours links instead of ugly boxes
	urlcolor     = blue, %Colour for external hyperlinks
	linkcolor    = blue, %Colour of internal links
	citecolor   = blue %Colour of citations
}

\usepackage[all]{xy}

\usepackage[shortlabels]{enumitem}

\newtheorem{thm}{Theorem}%[section]
\newtheorem {prop}{Proposition}
\newtheorem {lem}{Lemma}
\newtheorem {defi}{Definition}
\newtheorem {exa}{Example}

\newtheorem {cor}{Corollary}

\newtheorem {conj}{Conjecture}

\usepackage[capitalise]{cleveref}

\crefname{thm}{Thm.}{}
\crefname{prop}{Prop.}{}
\crefname{lem}{Lem.}{}
\crefname{cor}{Cor.}{}
\crefname{prob}{Prob.}{}
\crefname{figure}{Fig.}{}
\crefname{equation}{Eq.}{}
\crefname{exa}{Exa.}{}
\crefname{conj}{Conjecture}{}
\crefname{defi}{Def.}{}
\crefname{rem}{Remark}{}
\crefname{section}{Sec.}{}
\DeclareMathOperator\wgcd{wgcd \, }					%  weighted gcd 
\DeclareMathOperator\hgcd{hgcd \, }						% generalized  gcd
\DeclareMathOperator\hwgcd{h_{wgcd} \, }				% generalized weighted gcd
\DeclareMathOperator\wh{\mathcal S}   					%  weighted height
\DeclareMathOperator\lwh{\mathfrak s}   					%  logarithmic weighted height

\DeclareMathOperator\mult{mult}
\DeclareMathOperator{\lcm}{lcm}

\DeclareMathOperator\Pic{Pic}

\DeclareMathOperator\Supp{Supp}
\DeclareMathOperator\Cl{Cl}
\DeclareMathOperator\CaDiv{CaDiv}
\DeclareMathOperator\WeDiv{WeDiv}
\DeclareMathOperator\ord{ord}
\DeclareMathOperator\disc{Disc}
\DeclareMathOperator\sing{Sing}
\DeclareMathOperator\Exc{Exc}
\DeclareMathOperator\red{red}

\DeclareMathOperator\Proj{Proj}
\DeclareMathOperator\Spec{Spec}

\DeclareMathOperator\wnsc{wnsc}
\DeclareMathOperator\wnklt{wnklt}
\DeclareMathOperator\wnlc{wnlc}
\newcommand\A{\mathbb A}
\newcommand\N{\mathbb N}
\newcommand\Z{\mathbb Z}
\newcommand\Q{\mathbb Q}
\newcommand\R{\mathbb R}

\newcommand\X{\mathcal X}           
\newcommand\Ac{\mathcal A}

\newcommand\cC{\mathcal C}
\newcommand\M{{\mathcal M}}

\newcommand\cI{\mathcal I}

\newcommand\Nc{{\mathcal N}}
\newcommand\Zc{\mathcal Z}
\newcommand\Rc{\mathcal R}
\newcommand\U{\mathcal U}

\newcommand\Wc{\mathcal W}
\newcommand\cO{\mathcal O}
\newcommand\Y{\mathcal Y}
\newcommand\cL{\mathcal L}

\newcommand\w{\mathfrak q}
\newcommand\iso{\cong}
\newcommand{\kk}{{\bar k}}

\newcommand\hn{\mathfrak s}					%  Global weighted Weil  height
\newcommand\cH{\mathcal H} 

\newcommand\x{\mathbf x}
\newcommand\y{\mathbf y}
\newcommand\tX{\tilde \X}
\newcommand\tY{\tilde \Y}
\newcommand\tP{\tilde P}
\newcommand\il{\zeta}
\newcommand\tx{\tilde x}

\def\P{\mathbb P}
\def\a{\alpha}
\def\b{\beta}

\def\div{\mbox{div}}

\newcommand{\tony}[1]{{\color{Green} \sf $\spadesuit\spadesuit\spadesuit$ Tony: [#1]}}

\begin{document}

\title[Vojta's conjecture on weighted projective varieties]{Vojta's conjecture on weighted projective varieties}


\author{Sajad Salami}
\address{Institute of Mathematics,  Statistics, State University of Rio de Janeiro,  Rio de Janeiro, Brazil}
\email{sajad.salami@ime.uerj.br}

\author{Tony  Shaska}
%\address{Research Institute of Science and Technology (RISAT), Vlora, Albania}
\address{Department of Mathematics \& Statistics, Oakland University, Rochester Hills, MI}
%\email{shaska@oakland.edu}
\email{tanush@umich.edu}
%\email{shaska@risat.org}

\date{\today}

\begin{abstract}
We formulate Vojta's conjecture for smooth weighted projective varieties, weighted multiplier ideal sheaves, and weighted log pairs and prove that all three versions of the conjecture are equivalent. 
In the process,  we introduce  generalized weighted general common divisors   and express them   as  heights of  weighted projective spaces   blown-up  relative   to an  exceptional divisor.  
Furthermore, we  prove that assuming Vojta's conjecture for weighted projective varieties one can bound the $\log {\rm h_{wgcd}}$ for any subvariety of codimension $\geq 2$ and a finite set of places $S$. 
An analogue result is proved for weighted homogeneous polynomials with integer coefficients.
\end{abstract}

\maketitle

\setcounter{tocdepth}{3}

%\tableofcontents

\section{Introduction}
%^^^^^^^^^^^^^^^^^^^^^^^^^^^^^^^^^^*
The theory of      weighted local and global   heights   for   weighted projective varieties and closed subvarieties  was introduced in 
\cite{Salami-Shaska1}.  Weighted heights provide an efficient way of dealing with rational points on weighted varieties compared to  the classical  approach of considering their projective image under the Veronese map.   Some of the natural questions on this new approach is whether the main conjectures of arithmetic geometry can be formulated  for weighted varieties and weighted heights.  This paper accomplishes exactly that for one of the most important conjectures of arithmetic geometry, that of Vojta's conjecture. 

We generalize  the  generalized  logarithmic greatest common divisors      to   
generalized  logarithmic weighted greatest common divisors and  express them   as a weighted height of a blow-up relative to a weighted exceptional divisor.  Moreover, we are able to state Vojta's conjecture for smooth weighted projective varieties in terms of weighted heights. 

This paper is organized as follows:  In Sec. \ref{sec-2}, we give a quick summary of some of the basics of weighted projective spaces $\P_{\w, k}^n$, their weighted heights, and the concept of  rational points on weighted varieties.
In the last part of Sec. \ref{sec-2},
we provide an example of a weighted hypersurface having no rational points of weighted height $\leq 2$ in $\P_{(2,4,6,10), \Q}^4$.  This corresponds to the locus $\cL_2$ of genus two curves with extra involutions in weighted projective space in $\P_{(2,4,6,10), \Q}^4$; see \cite{sha-sha} for more details on this topic. 
The analogue statement in terms of projective varieties is that the  projective image $\phi_{60} (\cL_2)$, under the Veronese map $\phi_{60}$,  has no rational points of projective height $\leq 2^{60}$. 
This is to illustrate that the bound on the height  
changes so dramatically that it is well worth it to consider arithmetic properties of the weighted variety instead of its projective counterpart. 

In Sec. \ref{Sec-Voj}, we  review the 
Vojta's conjecture on algebraic points on projective varieties; see Conj. \ref{voj1},  We  state the conjecture using a correction term involving a multiplier ideal sheaf instead of using the normal crossing divisors; see Conj. \ref{voj2}, and summarize \cite{Yasuda2018} on Vojta's conjecture for log pairs. This makes it possible to drop the condition that the variety be smooth in the statement of the conjecture. Vojta's conjecture for log pairs is stated in Conj. \ref{Yasuda}.

In Sec. \ref{sect-5}, we  investigate whether it is possible to have analog statements for weighted varieties and weighted heights.   In  Conj. \ref{voj1wsp},  we state Vojta's conjecture for $\X$ a smooth weighted projective variety, $K_\X$ a canonical divisor, {$\Ac$} an ample divisor and $D$ a normal crossings divisor on $\X$, all defined over $\Q$.   An analogue of Conj. \ref{voj2}, for weighted projective varieties is stated in Conj. \ref{wpv-voj2}. The terminology and theory for weighted log pairs is developed in this section, so we are able to state Vojta's conjecture for weighted log pairs in Conj. \ref{voj6}.  Finally, in Cor. \ref{eq-conj}, we prove that      Conj. \ref{voj1wsp},    Conj. \ref{wpv-voj2},  and Conj. \ref{voj6} are equivalent.

In Sec. \ref{sect-6},   we   extend the concept of the \emph{generalized greatest common divisor}  
to that of  \emph{generalized weighted greatest common divisor}.  Furthermore, we prove that   generalized logarithmic weighted greatest common divisor   is equal to  weighted height of points on a blow-up of $\P_{\w, \Q}^n $ with respect to the exceptional divisor of the blow-up.  
%
%We prove that the generalized logarithmic weighted greatest common divisor $\log \hwgcd (\x) >0$  if and only if $\x \notin \sing (\P_{\w, \Q}^n)$  (cf. \cref{sing1}) and  
%analogues of   Theorems 1, 2, and 6 in \cite{MR2162351} for the weighted \gcds that are all subject to the validity of Vojta's conjecture for weighted  projective varieties. 

More precisely, we prove (cf. Thm. \ref{main0}) that  for  $\X $ a smooth weighted variety  in  $\P_{\w, \Q}^n$, $\Ac$ an ample divisor on $\X$,  $\Y \subset \X $ a smooth subvariety of codimension $r\geq 2$,   and  $-K_\X$  a normal crossing divisor whose support does not intersect $\Y$, assuming   Conj. \ref{voj1wsp},  for every finite set of places $S$ and every $0 < \varepsilon < r-1$ there is a proper closed subvariety   $ \Zc = \Zc (\varepsilon, \X, \Y, \Ac, S)  $ of $\X$, and  constants  
\[
C_\varepsilon = C_\varepsilon (\X, \Y, \Ac, S) \quad  and      \quad \delta_{\varepsilon} = \delta_{\varepsilon} (\X, \Y, \Ac),
\]
such that  for all  $P\in (\X \setminus \Zc)  (\Q)$ one has
\begin{equation}\label{meq10}
	\log \hwgcd (P; \Y) \leq \varepsilon \lwh_{\X, \Ac} (P) + \frac 1 {r-1+\delta_\varepsilon} \lwh^\prime_{\X, -K_\X, S} (P) + C_\varepsilon.
\end{equation}

Let $\w=(q_0, \cdots, q_n)$  be  a well-formed set of weights,    $\P_{\w, \Q}^n$ the weighted projective space over $\Q$ with weights $\w$, 
$m=\lcm(q_0, \cdots, q_n)$, and  $\Zc \subset \P_{\w, \Q}^n$  be a closed subvariety defined by $f_1, \cdots, f_t \in \Z_\w [x_0, \ldots , x_n]$,    such that $\Zc \cap  \sing(\P_{\w, \Q}^n )=\emptyset $,  with    codimension $r=n- \dim(\Zc) \geq 2$  in $\X$.
Assume that  $S$ is a   finite set of  primes, $\varepsilon >0$, and  the  Vojta's conjecture   holds for smooth weighted varieties (see Conj. \ref{voj1wsp}). Then,
there exists a nonzero weighted  polynomial $g \in \Z_\w[ x_0, \ldots , x_n]$     and a constant $\delta=\delta_{\varepsilon, \Zc} >0$,    such that every $\tilde \a= (\a_0, \cdots, \a_n) \in \Z^{n+1}$ with    $\wgcd(\a_0, \cdots, \a_n)=1$ satisfies either $g(\tilde \a)=0$ or
\begin{equation}
	\label{mm1}
	\gcd(f_1({\tilde \a}), \cdots, f_t({\tilde \a}))\leq  \max \, \left\lbrace  |\a_0|^{\frac{1}{q_0}}, \cdots,  |\a_n|^{\frac{1}{q_n}}\}\right\rbrace^\varepsilon  \cdot \left(  |\a_0  \cdots \a_n|'_S\right)^{\frac{1}{ m (r-1+ 		\delta  )}},
\end{equation}
where $|\cdot|'_S$ is the ``prime-to-$S$'' part of its origin (cf. Thm. \ref{main1}).

As a final remark, we mention that the Weighted heights provide a very natural approach of bounding the "size" of points on a weighted variety and the classical projective height is only a special case of weighted heights. 
It is still unclear whether this new approach of bounding the size of points on  weighted projective varieties via weighted heights,  instead of the traditional approach of bounding the size of points on projective varieties  will lead to any progress toward proving Vojta's conjecture.  That remains object of future research.

%********************************************************************************************
\section{Preliminaries on weighted projective varieties}\label{sec-2}
In   this section we provide some background on weighted projective varieties and their heights. For details the reader can consult   \cite{Salami-Shaska1} and its references. 
\subsection{Weighted projective varieties}
%
%We assume the reader is familiar with weighted projective varieties in the level covered in \cite{b-g-sh} and \cite{Salami-Shaska1}.   Let's recall some basic terminology.  

Given any  integer $n\geq 1$,  let $\w =(q_0, \dots, q_n)$ be a vector of positive integers.  
For any field $k$,   $\A^n_k$ and  $\P^n_k$ denote the affine and projective spaces  of dimension $n$ over $k,$ respectively.
Consider the polynomial ring $R=k_\w [x_0, \dots , x_n]$ where   $x_i$ has weight $  q_i$ for   $i=0,1,\cdots, n.$
%
\iffalse
% Every polynomial is a sum of monomials $x^d= \prod x_i^{d_i}$ with weight $\sum_{i=1}^n q_i d_j $.
%
%
%For  every $\lambda \in k^\ast$ and any  weighted homogeneous polynomial $f$ of degree $d$,   we have
%
%\[  f(\lambda^{q_0} x_0, \lambda^{q_1} x_1, \dots , \lambda^{q_n} x_n) = \lambda^d f(x_0, \dots , x_n).\]
%
For example,  a binary weighted form of   degree $d$, where   $w = (q_0, q_1)$ be  respectively the weights of $x_0$ and $x_1$, is  given by a polynomial as follows
%
\[
f(x_0, x_1) = \sum_{d_0, d_1} a_{d_0, d_1}  x_0^{d_0} x_1^{d_1},  \, \, \text{ such that }   \, \, d_0 q_0+ d_1 q_1 = d 
\]
%
and in decreasing powers of $x_0$ we have
%
\[ f(x_0, x_1)  =  a_{d/q_0, 0}  x_0^{d/q_0} + \dots + a_{d_0, d_1} x_0^{d_0} x_1^{d_1} + \dots+  a_{0, d/q_1} x_1^{d/q_1}\] 
%
By dividing   with $x_1^{d/q_1}$ and making a change of coordinates $X = x_0^{q_1} / x_1^{q_0}$ we get 
%
\begin{equation}
	f(x_0, x_1)  =    a_{d/q_0, 0} X^{d/q_0q_1} + \dots +  a_{d_0, d_1}  X^{d_0/q_1}+ \dots+  a_{0, d/q_1} = f(X)
\end{equation}
%
as noted in \cite{b-g-sh}.
% Notice that the condition $f(P)=0$ is well defined on $\P^n_{\w, k}$.  
%\end{exa}
\fi
%
We define   the  \textbf{weighted projective space} of dimension $n$ to be  the projective scheme 
$\Proj  \left(R \right)$, and denote it by  	$\P_{\w, k}^n$. 
Notice that if $\w=(1, \dots , 1)$, then 	$\P_{\w, k}^n$ is nothing but the usual projective space    $\P^n_k$. 
We simply write $\A^n$, $\P^n$, and $\P_{\w}^n$, for the affine, projective and weighted projective spaces of dimension $n$ over $\bar k$ an algebraically closed field containing $k$. 

We  say that $\P_{\w, k}^n $ is   \textbf{reduced} if   $\gcd (q_0,  \cdots , q_n)   = 1$.   It is called  \textbf{normalized} or \textbf{well-formed} if   
$ \gcd (q_0, \dots , \hat q_i, \dots , q_n)   = 1, $ for each $i=0, \dots , n.$
We notice that  any $\P^n_{\w, k}$  is thus isomorphic to a reduced and  well-formed one.  
To see the proof, one can see \cite[Prop. 3.3]{Salami-Shaska1} or \cite[Sec. 3C]{Beltrametti1986}.

A polynomial   $f \in R$   is called a \textbf{weighted homogeneous of  degree  $d$} if every monomial of $f$ has weight $d$, i.e. 
\[ f (x_0, \dots , x_n) = \sum_{i=1}^t a_i \prod_{j=0}^n x_j^{d_j}, \; \text{ for } \;  a_i \in k  \, \text{ and } \,  t \in \N, \]
and for all $ 0 \leq i \leq n$, we have  $\sum_{i=1}^n q_i d_j = d.$

By  a \textbf{weighted variety} we mean    an integral, separated subscheme $\X$ of finite type of $\P_{\w, k}^n$, 
equivalently,   if there are weighted homogeneous polynomials
$f_1, \cdots, f_t \in R$ such that $\X$ is isomorphic  to the $k$-scheme 
$\Proj \left( \frac{R}{\left\langle f_1, \cdots, f_t  \right\rangle }\right).$
For example,  a \textbf{weighetd hyperplane} in $\P_{\w, k}$  is defined as  
$\Proj \left( \frac{R}{\left\langle \ell \right\rangle }\right),$
where $\ell$ is a weighted homogeneous polynomial of degree $m=\lcm (q_0, \ldots , q_n)$ of the form 
\begin{equation}
	\ell(x_0, \dots, x_n)= a_0 x_0^{m/q_0}+ a_1 x_1^{m/q_1} + \cdots+ a_n x_n^{m/q_n},
\end{equation}
When     the  positive  integer $m$ is an arbitrary multiple of $\lcm (q_0, \ldots , q_n)$,   above equation defines a \textbf{weighetd diagonal hypersurface} in $\P_{\w, k}$ of degree $m$, see \cite{Go1996}.

%^^^^^^^^^^^^^^^^^^^^^^^^^^^^^^^^^^^
%\subsection{Veronese map}
%
%	For any     graded ring  $R=\bigoplus_{j \geq 0} R_{ j}$  and  integer $d \geq 1$, its \textbf{$d$-th truncated ring} is  
%	the subring $R^{[d]} \subseteq R$  defined by  $R^{[d]} = \bigoplus_{j \geq 0} R_{d j}.$
%
%	Clearly we have  the embedding $R^{[d]}  \hookrightarrow  R$, which is called the \textbf{$d$-th Veronese embedding},
%	implying   that  $\Proj (R^{[d]})  \iso  \Proj (R)$ by  \cite[Prop.~ 2.4.7]{Grothendieck1961}.
%	Based on this fact one can show that any weighted projective space  $\P^n_{\w, k}$ is isomorphic to  $\P^n_{\w',k}$, where  $\w'$ is a reduced tuple of weights.
%	If  $\P^n_{\w, k}$  is reduced and  $d_i= \gcd (q_0, \cdots, \hat q_i, \cdots, q_n)$ for $0\leq i \leq n$, 
%	then $\P^n_{\w, k} \iso  \P^n_{\w', k} $ with  $\w'=\left(\frac{q_0}{d_i}, \dots, \frac{q_{i-1}}{d_i}, q_i, \frac{q_{i+1}}{d_i}, \dots, \frac{q_n}{d_i}\right).$
%
%	Moreover, any $\P^n_{\w, k}$  is thus isomorphic to a reduced and  well-formed one.  
%
%To see the proof, one can see \cite[Prop. 3.3]{Salami-Shaska1} or \cite[Sec. 3C]{Beltrametti1986}.

%^^^^^^^^^^^^^^^^^^^^^^^^^^
\subsection{Rational points on weighted varieties in $\P_{\w}^n $}
For  $\w=(q_0, \dots , q_n)$,  we assume that $k$ contains a fixed primitive $q_i$-th root of unity $\xi_{q_i},$ which generates a cyclic subgroup  of $k$ denoted by  $\mu_{q_i}$  for every $i=0, \cdots, n$.
Let    
$G_\w\colonequals \mu_{q_0} \times \cdots \times \mu_{q_n}$,  which is a finite group of order $|G_\w|= \prod_{i=0}^{n} q_i $.  Then, there is an action of $G_\w$ on $\P^n$  given by 
\begin{equation}\label{equivalence2}
	(\xi_0,  \cdots, \xi_n) \bullet [y_0: \dots : y_n]   = \left[ \xi_0  y_0: \dots : \xi_n  y_n   \right].
\end{equation}
The weighted projective space $\P_{\w}^n$ is isomorphic  to the quotient space  of this action, i.e.,
say $\P_{\w}^n \cong \P^n /G_{\w} $    as described below.
%Note that $G_\w \iso \mu_m$ if and only if  $m = \lcm (q_0, \dots , q_n)$, that is,  all of $q_i$'s are pairwise coprime. 
%In this case, action of $G_m$    on $\P_k^n$ can be expressed as  
%
%\begin{equation}\label{equivalence3}
%	\xi^\alpha \cdot  [x_0: \dots : x_n]   =
%	\left[ \xi^{\alpha /q_0} x_0: \dots : \xi^{\alpha /q_n} x_n   \right],
%\end{equation}
%
%for $0 \leq \alpha \leq m-1$, where $\xi \in G_m$ is a  $m$-th root of unity.
Let  $\textbf{0}=(0, \cdots, 0 )$ be the origin of $ \A^{n+1}$, and consider the  morphism  
$	\pi_0: \A^{n+1} \setminus \{{\bf 0}\}     \longrightarrow  \A^{n+1} \setminus \{{\bf 0}\} $ defined by
\begin{equation} 
	\pi_0(y_0, \cdots, y_n)  = \left(x_0, \cdots, x_n \right) :=\left( y_0^{q_0}, \dots , y_n^{q_n}   \right).
\end{equation}
It	induces  the following diagram
\begin{equation}  \label{diag1}
	\xymatrixcolsep{2pc}
	\xymatrix{ 
		\A^{n+1} \setminus \{{\bf 0}\} \ar[rr]^{\pi_0} \ar[d]^{p}  &  & \A^{n+1} \setminus \{{\bf 0}\} \ar[d]^{p_\w} \\
		\P^n \ar[rr]^{\pi_{\w}}   \ar[dr]^{\bar{\pi}_\w}  &   &  \P_{\w}^n \\
		& \P^n/ G_\w  \ar[ur]_{\iso } & }
\end{equation}
where   $p:  \A^{n+1} \setminus \{{\bf 0}\}   \to  \P^n$ is the canonical  map
$(y_0, \cdots, y_n) \mapsto [y_0: \cdots: y_n] $, and 	
$			p_\w: \A^{n+1} \setminus \{{\bf 0}\}         \longrightarrow  \P_{\w}^n  $ is given by
\begin{equation*} 
	p_\w \left( x_0, \cdots, x_n	\right) =	p_\w \left( y_0^{q_0}, \dots , y_n^{q_n}  \right)   = \left[ y_0^{q_0} : \dots : y_n^{q_n}   \right].
\end{equation*}	 
Hence, $\bar{\pi}_\w:  \P^n \longrightarrow   \P_{\w}^n/ G_\w $ is the induced quotient map of the action, which commutes with  		
the map $	\pi_\w:  \P^n        \longrightarrow    \P_{\w}^n  $ given by 
\begin{equation} 
	\pi_\w \left([y_0: \cdots: y_n] \right)    = [x_0: \cdots: x_n]:=\left[ y_0^{q_0} : \dots : y_n^{q_n}   \right].
\end{equation}	 	
The morphism $\pi_{\w} $ is surjective, finite, and  its   fibers  are orbits of the action of $G_\w$ on $\P^n$.
See   \cite[Chap. V, Props. 1.3 and  1.8]{Grothendieck1961} for more details.
Furthermore, 	$\P_{\w}^n$ 	 is    an irreducible,  normal and Cohen-Macaulay  variety having  only    cyclic quotients singularities, see \cite{dolga}.

Conversely, there is a map $	\phi_m : \P^n_{\w}      \longrightarrow  \P^n$ which is defined by
\begin{equation}\label{veronese}
	\phi_m([x_0: \cdots : x_n])  = [y_0: \cdots : y_n]:=[x_0^{m/q_0} :  x_1^{m/q_1} :  \cdots : x_n^{m/q_n}].
\end{equation}
We notice that  $\phi_m$ and $\pi_{\w}$ are non-canonical maps between $\P^n_{\w} $ and  $\P^n$. 	
If $\w$ is reduced and all of   $m/q_i$ are co-prime, where $  m= \lcm \left( q_0, \cdots ,  q_i\right),$  then $\P^n_{\w}$ is isomorphic to $\P^n$ by   via the map $\phi_m$ and its inverse  is  $\pi_{\w}$, see \cite[Sec. 3C]{Beltrametti1986}.
%

% 
%	We  call  the isomorphism  $\phi_m$ given in \eqref{veronese}   the \textbf{Veronese map}; see \cite[Prop. 3]{Salami-Shaska1} for details. 

As in the usual projective spaces, a point $\x= [x_0: \cdots : x_n] \in \P_{\w}^n $ is   $k$-rational  if and only if $p_\w^{-1}(\x)$ is invariant under the Galois group $\operatorname{Gal}(\bar{k} / k)$. We denote  the subset of rational points on $\P_{\w}^n $ by $\P_{\w}^n(k)$. The orbit of any $x=\left(x_{0}, \ldots, x_{n}\right) \in  	\A^{n+1} \setminus \{{\bf 0}\} $ with image $p_\w(x)=\x$ is the rational curve
\[
\cC^*_x=p_\w^{-1}(\x)= \left\{\left(\lambda^{q_{0}} x_{0}, \ldots, \lambda^{q_{n}} x_{n}\right): \lambda \in \bar{k}^{\times}\right\} \subset \A^{n+1} \setminus \{{\bf 0}\}.
\]
Moreover,   
one has $p_\w(x) \in \P_{\w}^n(k) $ if and only if  
$\cC^*_x  \cap  \left( 	\A^{n+1} \setminus \{{\bf 0}\}\right)  (k) \neq \emptyset$.
In other words, the map 
\begin{equation}
	p_\w:  \left( 	\A^{n+1} \setminus \{{\bf 0}\}\right)  (k)   \longrightarrow \P_{\w}^n(k)
\end{equation}
is surjective and   induces a bijection 
\[
\left( 	\A^{n+1} \setminus \{{\bf 0}\}\right)  (k) /\sim \ \  \xrightarrow{\sim} \P_{\w}^n(k),
\]
where $ \sim$  denotes an equivalence relation on $\left( \A^{n+1} \setminus \{{\bf 0}\}\right)  (k)$ defined as follows:
for two points $x =(x_0 , \dots, x_r) $ and $x' =( x'_0 , \dots, x'_n)$ in $\left( \A^{n+1} \setminus \{{\bf 0}\}\right)  (k)$, we write
$x\sim x'$ if there exists some $ \lambda \in   {\bar k}^\ast$ such that $(x_0 , \dots, x_r) =(\lambda^{q_0} x'_0 , \dots, \lambda^{q_n} x'_n)$.
These can be proved by Hilbert's Theorem 90, as in the usual projective  case,
for more details one can see  \cite[Lemma 1.2]{Go1996}.

For any  weighted subvariety     $\X \subseteq \P_{\w}^n $, 
the  \textbf{punctured affine cone}  over $\X$ is  defined  as $\cC^*_\X=p_\w^{-1} (\X)$. 
The \textbf{affine cone}   $\cC_\X$ over $\X$ is the closure of $\cC^*_\X$ in $\A^{n+1} $.  
The origin point $\textbf{0}$ refers to the vertex of  $\cC_\X$.
%We note that ${\bar k}^\ast$ acts on the  punctured affine cone $\cC^*_\X=p_\w^{-1} (\X)$  to result  $\X= \cC^*_\X / {\bar k}^\ast$.

%A  weighted subvariety   $\X$  of  $\P_{\w}^n $ is called \textbf{quasi-smooth}  of dimension $m$   if its affine cone  $\cC_\X$ is smooth   variety of dimension $m+1$ outside its vertex.     The singularities of a quasi-smooth variety $\X$  are due to the ${\bar k}^\ast$-action  and hence are cyclic quotients singularities.
%Furthermore,  if $\X \subset \P^n_{\w}$ is subvariety such that $\X \cap \sing (\P^n_{\w})= \emptyset$, then  $\X$ is non-singular if and only if $\X$ is quasi-smooth.

A point $\x \in \X$ is a $k$-rational point on $\X$ if $\iota(\x)=[x_0: \dots: x_n] \in \P_{\w}^n(k),$ where
$\iota: \X \rightarrow  \P_{\w}^n$ is the natural embedding.  We denote by $\X(k)$ the set of all $k$-rational points on $\X$.
Thus, the following diagrams commute 

\begin{equation}  \label{diag2}
	\xymatrixcolsep{2pc}
	\xymatrix{ 
		\cC^*_{\X} \ar[rr] \ar[d]^{p_\w}  &  & \A^{n+1} \setminus \{{\bf 0}\} \ar[d]^{p_\w} \\
		\X \ar[rr]^\iota    &   &  \P_{\w}^n  } 
	\quad \quad 
	\xymatrixcolsep{2pc}
	\xymatrix{ 
		\cC^*_{\X} (k) \ar[rr] \ar[d]^{p_\w}  &  & \left( \A^{n+1} \setminus \{{\bf 0}\}\right)  (k) \ar[d]^{p_\w} \\
		\X(k) \ar[rr]^\iota     &   &  \P_{\w}^n(k)  }		
\end{equation}

%**************************
\subsection{Weighted heights on  $\P_\w^n$ }
Let  $k$ be an algebraic number field as in the previous section
and $\kk$ be an algebraically closed field containing $k$. 
%We denote by  $\cO_k$   the  ring of algebraic  integers in $k$. 
%Let   $\X$ be a variety over $k$, i.e.   an integral separated scheme of finite type over $\Spec(k)$ and  $\cO_\X$ the ring sheaf of regular  functions on $\X$.
% We will use $\X$ to mean $\X(\kk)$ and   $\X(k)$  for the set of $k$-rational points on $\X$. 
Denote by $M_k$ the set of all places of $k$, i.e. the equivalent  classes of  absolute values on $k$. It is a disjoint union of $M_k^0$, the set of all non-archimedian places,  and  $M_k^\infty$,   the set of all Archimedean places of $k$.   More precisely, if $\nu \in M_k^0$, then 
%$\nu=\nu_\p$ for some prime ideal $\p \subset \cO_k$  over a prime number  $p$  such that  
$\nu|_\Q $ is the $p$-adic absolute value  $|\cdot|_p$ on $\Q$.  If $\nu \in M_k^\infty$, then $\nu=\nu_\infty$ and $\nu_\infty|_\Q$ is the usual absolute value $|\cdot|_\infty$ on $\Q$.

Given  any $\x  \in \P_\w^n (k)$, the \textbf{multiplicative weighted  height} over $k$  is defined as  
\begin{equation}\label{def:height}
	\wh_k ( \x ) \colonequals \prod_{\nu \in M_k} \max   \left\{   |x_0|_\nu^{\frac{1}{q_0}} , \dots, |x_n|_\nu^{\frac{1}{q_n}}  \right\}
\end{equation}
and its  \textbf{logarithmic  weighted  height} (over $k$) as 
\begin{equation}\label{log-height}
	\lwh_k(\x) \colonequals \log \wh_k (\x)=    \sum_{\nu \in M_k}   \max_{0 \leq j \leq n}\left\{\frac{1}{q_j} \cdot  \log  |x_j|_\nu \right\}.
\end{equation}

In \cite[Prop.~ 1]{b-g-sh} it is shown that height functions $\wh_k ( \x )$ and hence $\lwh_k(\x)$ are independent of the choice of   coordinates of the point $\x$. Moreover, in   \cite[Prop.~ 5-ii]{b-g-sh}, it is proved that for any finite extension $K|k$ we have 
\[
\wh_k ( \x )^{[K:k]}=\wh_K ( \x ), \; \text{and hence} \;  {[K:k]} \cdot \lwh_k (\x) = \lwh_K (\x).
\]
The \textbf{absolute weighted height} on $\P_\w^n(\kk)$ is defined as
\begin{equation}\label{wgh}
	\begin{split}
		\wh  : \P_\w^n(\kk)  &    \to [0, \infty],  \\
		\x    &   \mapsto    \wh (\x) \colonequals   \wh_K(\x)^{1/[K:k]},
	\end{split}
\end{equation}
and the \textbf{absolute logarithmic weighted height}   on $\P_\w^n(\kk)$ is given by
\begin{equation}\label{wgh1}
	\begin{split}
		\lwh  : \P_\w^n(\kk)    &   \to [0, \infty ],   \\
		\x   \mapsto      &  \lwh (\x)\colonequals \frac{1}{[K:k]} \log   \wh_K(\x),
	\end{split}
\end{equation}
for which $K\subset \kk$ is a finite extension of $k$ containing   $k(\x),$ 
the
% \textbf{field of   normalization}  of    $\x$ is defined   by  
%$K\colonequals k(\hwgcd (\x))$
%and
\textbf{field of   definition }   of    $\x$  defined by
\[
k(\x)\colonequals k \left( \frac{x_0^{1/q_0}}{x_i^{1/q_i}} , \cdots, 1, \cdots,     
\frac{x_n^{1/q_n }}{x_i^{1/q_i }} \right), 
\]
for some $x_i\neq 0$.
Notice that both of these height functions are independent of the choice of the field $K$; see \cite{b-g-sh}. For simplicity, we call   $\lwh (\x)$  the  \textbf{global  weighted height} on $\P_\w^n (\kk)$.	For all $\x \in  \P_\w^n(\kk) $, we have
\begin{equation}
	%\begin{split}
	\wh ( \x )   = H  \left( {\phi_m} (\x) \right)^{\frac{1}{m}}, \; \text{ and } \;     \lwh(\x)  = \frac{1}{m} \cdot h \left( {\phi_m} (\x) \right),
	%\end{split}
\end{equation}
where   $\phi_m$ is as in \eqref{veronese},  $H (\cdot)$,  $h (\cdot)$ are multiplicative and logarithmic heights $\P^n$, and $\wh (\cdot)$, $\lwh (\cdot)$ are defined as in  \eqref{wgh} and \eqref{wgh1}.

Since the weighted projective spaces admit natural stack structures,  the above definition of weighted heights  
is    compatible with   the stable height  introduced   by  Darda  in\cite{Dar21},  Ellenberg–Satriano–Zureick-Brown in \cite{ESZB-2023}, and Darda-Yasuda in \cite{Dar-Yas24} to explore the distributions of rational points on algebraic stacks as part of their efforts to unify the Batyrev–Manin conjecture and Malle’s conjecture. To achieve this, they introduced new height functions on stacks, drawing from distinct theoretical frameworks. Specifically, Darda \cite{Dar21} proposed the concept of quasi-toric heights on weighted projective stacks. In \cite{ESZB-2023},  the authors developed a height function linked to vector bundle on a proper algebraic stack with finite diagonal over a global field. In the current work \cite{Dar-Yas24},  Darda and Yasuda  introduced  a new height function on rational points of a  Deligne–Mumford  stack over a number field,  which  generalizes  those in   \cite{Dar21, ESZB-2023}.

We end up this section by  providing an example of a weighted hypersurface and its rational points of small weighted height illustrating that it would be difficult to obtain this result by embedding this to a projective space counting rational points on that projective space.

\begin{exa}
	Consider $\w=(2,4,6,10)$ and a weighted  hypersurface in $\P^3_{\w, \Q}$ given by the following equation
	\begin{small}
		\[
		\begin{split}
			&	41472 w y^5+159 y^6 x^3-236196 w^2 x^5-80 y^7 x+104976000 w^2 x^2 z-1728 y^5 x^2 z+6048 y^4 x z^2\\
			&	 -9331200 w y^2 z^2-2099520000 w^2 y z+12 x^6 y^3 z-54 x^5 y^2 z^2+108 x^4 y z^3+1332 x^4 y^4 z\\
			&	-8910 x^3 y^3 z^2+29376 x^2 y^2 z^3-47952 x y z^4-x^7 y^4-81 x^3 z^4-78 x^5 y^5+384 y^6 z-6912 y^3 z^3\\
			&	+507384000 w^2 y^2 x-19245600 w^2 y x^3-592272 w y^4 x^2+77436 w y^3 x^4+4743360 w y^3 x z\\
			&	-870912 w y^2 x^3 z+3090960 w y x^2 z^2-5832 w x^5 y z-125971200000 w^3+31104 z^5+972 w x^6 y^2\\
			&	+8748 w x^4 z^2-3499200 w x z^3   = 0
		\end{split}
		\]
	\end{small}
	This is the locus $\cL_2$  of genus two curves 	with extra involution where $(x, y, z, w) = (J_2, J_4, J_6, J_{10})$, see  \cite{sha-sha}. 
	The map  $\phi_{60} : \P^3_{\w, \Q} \to \P_\Q^3 $ is given by  
	\[
	\x=[x: y: z: w] \to [ x^{30}: y^{15}: z^{10}: w^6].
	\]
	There is a one to one correspondence between the rational points $\cL_2 (\Q) \subset \P_{\w} (\Q)$ and the rational points of $\phi_{60} (\cL_2) \subset \P^3_\Q$.  An inverse of this map is given in \cite{univ}. 
	
	In \cite{sha-sha}, it  is proved that the above weighted hypersurface has no rational points of weighted height $\wh  \leq 2$.  
	To prove this result via the projective varieties 	we have to check the rational points of $\phi_{60} (\cL_2) \subset \P^3_\Q$. Since, 
	$H (\phi_m (\x) ) = \wh (\x)^m$, this is equivalent to checking for rational points of projective height   $H (\phi_m (\x) )  \leq 2^{60}$.  
	
	The example above illustrates yet again how efficient computations in the weighted projective spaces are compared to the projective space.
	\qed
\end{exa}

%^^^^^^^^^^^^^^^^^^^^^^^^^^^^^^^^^^^*
\section{Vojta's conjecture for projective varieties}\label{Sec-Voj}

In this section, we recall the Vojta's famous conjecture on   algebraic points of  projective varieties, as well as its variants   using a correction term involving a multiplier ideal sheaf  and its version for log pairs of projective varieties.
The interested readers can find more in \cite{vojta-book, Vojta2011, Vojta2012, Yasuda2018}.

\subsection{Vojta's conjecture for smooth  varieties}
Here we give a brief overview of Vojta's conjecture over projective varieties defined over  $\Q$.
For any finite extension $L$ of    $\Q$,    we define the \textbf{logarithmic discriminant} $d (L)$ by
%
%\[ d_k (L) \colonequals   \frac   1 {[L:k]}    \log  |\disc(L)| -\log |\disc (k )|,\]
%
%
\[ d (L) \colonequals   \frac   1 {[L: \Q]}    \log  |\disc(L)|.\]
%
%where $\disc(\cdot)$ denotes the absolute discriminant.  
%
Given  a variety $\X$ over $\Q$  and a point  $\x \in \X$, we define its \textbf{logarithmic discriminant} by 
$d (\x)  \colonequals d(\Q(\x))$.  

Recall that a Cartier divisor $D$ on a smooth projective variety $\X$ is a \textbf{normal crossing divisor} if at every point in the support of $D$ there are local coordinates $z_0, z_1, \ldots, z_n$ such that $D$ is given locally by an equation of the form  $z_0 z_1 \ldots  z_n=0$.   Furthermore, the \textbf{canonical divisor} of $\X$ is a divisor $K_\X$ such that $\cO_\X (K_\X) = \omega_\X$, where 
$\omega_\X$ is the canonical sheaf of regular forms on $\X$.  Paul Vojta made a conjecture on  algebraic points on projective varieties     as follows:

\begin{conj}\label{voj1}
	Let $\X$ be a smooth projective variety defined over   $\Q$,   $K_\X$ a canonical divisor, $\Ac$ an ample divisor,  and $D$ a normal crossings divisor on $\X$.  Furthermore, let $S$ be a finite subset of places of $\Q$ containing $M_\Q^\infty$.  Then, given any real constant  $\varepsilon >0$ and any positive integer $r$, there exists  a proper Zariski-closed subset  $Z$ of $\X$, depending only on $  \X,  D,   \Ac,  \varepsilon,$ and $ r$, such that 
	\[
	h_{K_\X}(\x) + \sum_{\nu \in S}^{ } \lambda_D(\x, \nu) \leq  \varepsilon h_\Ac(\x) + d(\x) + O(1),
	\]
	for all $\x \in (\X \setminus Z) (\bar \Q)$ with $[\Q(\x):\Q]\leq r$.
\end{conj}
The case $r=1$ of the above conjecture  is   known in the literature
as Vojta's conjecture for the rational points of  algebraic varieties.

\subsection{Vojta's conjecture using  a multiplier ideal sheaf}

In \cite{Vojta2012}, Vojta restated his  conjecture using a correction term involving a multiplier ideal sheaf  instead of the normal crossing divisors as follows.

Let $\cI$ be a nonzero sheaf of ideals on a projective variety $\X$ and $c$ some positive  real constant. 
Let $f: \X' \rightarrow \X$  be a proper birational morphism such that $\X'$ is a smooth variety and  $f^* \cI = \cO_{\X'} (-E)$, for some normal crossing  divisor $E$ on $\X'$. Denote by $\Rc_{\X'/\X}$ the ramification divisor  of $\X'$ over $\X$ and  define the \textbf{multiplier ideal sheaf}  $\cI_c$ and  $\cI_c^-$ associated to $\cI$ and a rational number $c$ as 
\[ \cI_c\colonequals  f_* \cO_{\X'} (\Rc_{\X'/\X}  - \lfloor c E \rfloor),\ \text{and}\  \cI_c^- \colonequals\lim_{\varepsilon \rightarrow 0^+} \cI_{c-\varepsilon},\]
where $\lfloor c E \rfloor$  denotes the round down of the divisor $c E$.
We will denote $\cI_1$ and  $\cI_1^-$ by  $\cI$ and  $\cI^-$ respectively. 

\begin{conj}\label{voj2}
	Let $\X$ be a smooth projective variety defined over $\Q$,
	$K_\X$ a canonical divisor,  $\Ac$ an ample divisor,  and $\cI$  a nonzero ideal sheaf $\X$.
	Let $S$ be a finite subset of places of $\Q$  containing $M_\Q^\infty$. 
	Then, given any real constant  $\varepsilon >0$ and positive integer $r$, there exists a proper Zariski-closed subset  $Z$ of $\X$, depending only on $ \X,  \cI,   \Ac,  \varepsilon,  r$, such that 
	\[
	h_{K_\X}(\x) + \sum_{\nu \in S}^{ } \lambda_\cI(\x, \nu) -  \sum_{\nu \in S}^{ } \lambda_{\cI^-}(\x, \nu) \leq  \varepsilon h_\Ac(\x) + d(\x) + O(1).
	\]
	for all $\x \in (\X \setminus Z) (\bar \Q)$ with $[\Q(\x):\Q]\leq r$.
\end{conj}

%^^^^^^^^^^^^^^^^^^^^^^^^^^^^^^^^^^^*
\subsection{Vojta's conjecture for log pairs over $\Q$}\label{log-pairs}
Since   Conj. \ref{voj1}, and Conj. \ref{voj2},  do not deal with  singular varieties in \cite{Yasuda2018} Yasuda formulated a generalization of it  in terms of log pairs and variants of multiplier ideals.  In order to state his generalization, first we need to recall some   terminology. The reader can refer to  \cite{Kollar2013, Lazarsfeld2004c, Lazarsfeld2004b}, or \cite{Yasuda2018} for more details.

Let $\X$ be a variety defined over $\Q$. Then,  $\X$ is said to be \textbf{$\Q$-Gorenstein} if it  is  Gorenstein in codimension one,   satisfies Serre's condition $S_2$,   and a canonical divisor $K_\X$ is $\Q$-Cartier. For example, if $\X$ is normal, then the first two conditions are true automatically and hence a canonical divisor exists unequally up to linear equivalence  and  is Cartier in  codimension one.

A \textbf{$\Q$-subscheme} of  $\X$ is a formal linear combination $\Y=\sum_{i=1}^m c_i \cdot \Y_i$ of proper closed subschemes $\Y_i \subset \X$ with  all $c_i\in \Q$. The support  of such $\Y$ is defined to be the closed subset $\displaystyle \cup_{c_i \neq 0} \Y_i$, and  it is called   \textbf{effective $\Q$-subscheme} if $c_i \geq 0$ for every $i$. By a \textbf{log pair}, we mean a pair $(\X, \Y)$ of a  $\Q$-Gorenstein variety and an effective   $\Q$-subscheme $\Y$ of $\X$. For example, if $\X$ is a normal $\Q$-Gorenstein and $D$ is an effective $\Q$-divisor, then $(\X, D) $ is a log pair.

A \textbf{resolution} of $\X$ over $\Q$  is a projective birational morphism $f: \X' \rightarrow \X$ such that  
$\X'$ is a smooth projective variety  over $\Q$.    
By a \textbf{log resolution} of  a log pair $(\X, \Y)$  with $\Y=\sum_{i=1}^m c_i \cdot \Y_i$, we mean a resolution   $f: \X' \rightarrow \X$ of $\X$ such that the set-theoretic inverse image $f^{-1} (\Y_i)$  is a Cartier divisor on $\X'$, and  the union of  exceptional divisor $\Exc(f) $  of $f$ with all $f^{-1} (\Y_i)_{\red}$   is a simple normal crossing divisor of $\X'$.
For a log resolution  $f: \X' \rightarrow \X$  of a log pair $(\X, \Y)$, the \textbf{relative canonical divisor}  of $\X'$ over  $(\X, \Y)$ is defined to be the  $\Q$-Weil divisor 
\[
K_{\X'/(\X, \Y) }=K_{\X'/\X } - f^* \Y, 
\]
where  $K_{\X'/\X } $ is the  relative canonical divisor  of $\X'$ over $\X$, and $f^* \Y$ is the pull-back of $\Y$ by $f$ over $\X'$.

For a log pair $(\X, \Y)$ with a log resolution $f: \X' \rightarrow \X$, we define $\cI(\X, \Y)$ a \textbf{variant of multiplier sheaf} as
\[
\cI(\X, \Y)\colonequals f_* \cO_{ \X'} (\lceil K_{\X'/(\X, \Y) } \rceil)
\]
if $\X$ is a normal variety, otherwise, we let 
\[
\cI(\X, \Y)\colonequals \bar{f}_* \cO_{ \X'} (\lceil K_{\X'/(\X, \Y) } \rceil),
\]
where   $\lceil \cdot \rceil $  denotes the round up of the $\Q$-divisors, 
and 
$\bar{f}_*   \cO_{ \X'}(E)$ denotes the largest  ideal sheaf  in $\cO_{ \X}$ for which its pull-back by $f$ is contained in  $\cO_{ \X'}(E)$ as an $\cO_{\X'}$-submodule of   (constant) function field sheaf $\M_{\X'}$. 
Moreover, there exist a constant $\varepsilon_0>0$ such that for every rational number 
$0< \varepsilon \leq \varepsilon_0$, one has  $\cI(\X, (1-\varepsilon) \Y)=\cI(\X, (1-\varepsilon_0) \Y)$.
Based on this fact, we set 
\[\cI^-(\X, \Y)\colonequals\cI(\X, (1-\varepsilon) \Y), \  (0< \varepsilon \ll 1).\]
We also define another ideal sheaf as
\[
\cH(\X, \Y)\colonequals\bar{f}_* \cO_{ \X'} (\lfloor K_{\X'/(\X, \Y) } \rfloor), 
\]
where $\bar{f}_* $  is as above.
The definitions of $\cI(\X, \Y)$, and hence $\cI^-(\X, \Y)$, as well as $\cH(\X, \Y)$ are  independent of the choice of a log resolution   by  \cite[Lem. 3.1]{Yasuda2018} and  \cite[Prop. 3.4]{Yasuda2018} respectively.

For a $\Q$-Gorenstein projective variety  $\X$ with a canonical divisor $K_\X$, we can define  a global height
function $h_{K_\X}$ up to addition of a bounded function. For a log pair $(\X, \Y)$ of a $\Q$-Gorenstein $\X$, 
we define  the logarithmic height function  associated to the subscheme $K_{(\X, \Y)}:=K_{\X}+ \Y$ as  
\begin{equation}
	h_{K_{(\X, \Y)}}=h_{K_\X}+ h_\Y,
\end{equation}
where $h_\Y$ is the logarithmic  height function associated to  $\Y$ or its  ideal sheaf.
Next is   Yasuda's generalization of Vojta's conjecture for algebraic points.

\begin{conj}\label{Yasuda}
	Let  $(\X, \Y)$ a log pair with projective $\X$, $\Y$ a closed subscheme with ideal sheaf $\cI=\cI(\Y)$, $K_\X$ a canonical divisor and  $\Ac$ an ample divisor on $\X$  all defined over   $\Q$.
	Let $S$ be a finite subset of places  of $\Q$  containing $M_\Q^\infty$.  Then, given any real constant  $\varepsilon >0$ and positive integer $r$, there exists  a proper Zariski-closed subset  $Z$ of $\X$, depending only on $  \X,  \cI,   \Ac,  \varepsilon,  r$, such that 
	\[ 
	h_{K_{(\X, \Y)}} - \sum_{\nu \in S}^{ } \lambda_\cH(\x, \nu) -  \sum_{\nu \in S}^{ } \lambda_{\cI^-}(\x, \nu) \leq  \varepsilon h_\Ac(\x) + d(\x) + O(1).
	\]
	for all $\x \in (\X \setminus Z) (\bar \Q)$ with $[\Q(\x):\Q]\leq r$, where $\cH=\cH(\X, \Y)$ and $\cI^-=\cI^-(\X, \Y).$
\end{conj}

One can see that the above conjecture holds for a log pair  $(\X, \Y)$  and a log resolution  $f: \X'\rightarrow \X$,  if  Conj. \ref{voj1}, holds  for $\X'$ and  the reduced simple normal crossing divisor
\[ \lceil K_{\X'/(\X, \Y)}  +\varepsilon f^* \Y \rceil  - \lfloor  K_{\X'/(\X, \Y)} \rfloor,\]
for $0 < \varepsilon \ll 1$, where    $\lfloor \cdot \rfloor$  denotes the round down of the divisor.
Moreover, the Conjecture \ref{voj1}, Conjecture \ref{voj2},  and Conjecture \ref{Yasuda} are equivalent; see \cite[Prop. 5.4 and Rem. 5.5]{Yasuda2018} for a proof.

%^^^^^^^^^^^^^^^^^^^^^^^^*
%\newpage
\section{Vojta's conjecture for weighted projective varieties} \label{sect-5}
In this section, we will state the analogues of Vojta's conjectures given in   previous section for weighted projective varieties.
In order to do that we review  and fix some terminologies on weighted varieties.

\subsection{Local and global weighted heights}

Let   $\X$ be a weighted  variety over $k$, i.e.   an integral separated scheme of finite type over $\Spec(k)$.
In \cite[Prop. 4.3]{Salami-Shaska1},    It is proved  that any line bundle $\cL$ on   $\X$ admits a locally bounded weighted $M$-metric  $\| \cdot \|=(\|\cdot \|_u)$. One can see   \cite[Sec. 4 ]{Salami-Shaska1}, for more details on the subject. 
Given any   weighted Cartier divisor $D=\{ (U_i, f_i)\}$  on $\X$,  we let  $\cL_D=\cO_\X(D)$ be  the line bundle of regular functions on $D$, which can be constructed by gluing 
\[
\cO_\X(D)|_{U_i}= f_i^{-1} \cO_\X( U_i)
\]
and $1$ becomes a canonical invertible meromorphic section of $\cL_D$, which is denoted by $g_D$.
Thus,  we can equip
$\cL_D$   with a weighted  locally bounded $M$-metric $\| \cdot \|$,  determined    by the max-min method  in   proof of \cite[Prop. 4.3]{Salami-Shaska1}, and denote it  by $\widehat{D}=\left( \cL_D, \| \cdot \|\right).$

For each  $\nu \in M_k$, the \textbf{local weighted height  $\il_{\widehat{D}}(-,\nu)$ on $\X$  with respect to} $\widehat{D}$ on weighted variety  $\X$ is defined as
\[
\il_{\widehat{D}}(\x, \nu)=- \log \| g_D(\x)  \|_v,
\]  
for  $\x      \in \X \backslash  \Supp (D)$, 
where   $v\in M$ such that $\nu=v|_k$.   The main properties of  local weighted heights are proved in \cite[Thm. 4.4]{Salami-Shaska1}.
The \textbf{global weighted height $\hn_{\widehat{\cL}}(\x)$ } on $\X$   with respect to   $\widehat{\cL}$  is defined by 
\[
\hn_{\widehat{\cL}}(\x)\colonequals \sum_{u\in M_K}^{} \il_{\widehat{\cL_g}} (\x, u), 
\]
where   $\il_{\widehat{\cL_g}} (\x, u) = -\log \| g(\x)\|_u$, 
and its   properties  are described in \cite[Thm. 4.5]{Salami-Shaska1}.

Fix   a weighted  variety $\X $ in $\P_{\w, \kk}^n$ defined over $k$, i.e., an integral and separated subscheme of finite type.
Given any closed subscheme $\Y$ of $\X$ over $k$, we let $\cI_{\Y}$ denotes the corresponding 
sheaf of ideals generated by  $f_1,\cdots, f_r \in k_\w[x_0,  \cdots, x_n]$ with $\deg(f_j)=d_j$ for $j=1,\cdots r$.
Letting $D_j \colonequals \div(f_j)$ for $j=1, \cdots, r$ and  $\nu \in M_k$, 
we define    
\[
\il_\Y( \cdot , \nu): (\X \backslash  \Y) (k) \rightarrow \R,
\]
the \textbf{local weighted height on $\X$ associated to $\Y$}, by 
\begin{equation}\label{gwh1}
	\il_\Y(\x, \nu) \colonequals \min_{1\leq j \leq r}^{}\{ \il_{\widehat{D_j}}(\x, \nu) \} \\
	=\min_{1\leq j \leq r}^{} \left\lbrace   - \log \frac{|f_j(\x)|_\nu}{ \max_{i}^{} \left| x_i\right|_\nu^{\frac{ d_j}{q_i}}}\right\rbrace.
\end{equation}
By convention, we define $\il_\Y(\x, \nu)=\infty$  for $\x \in \Y (k)$.
%One can show that this is unique up to a weighted $M_k$-bounded function by a similar argument for the projective varieties.

The \textbf{global weighted height on $\X$ associated  to $\Y$}, can be defined  up to a bounded function by summing all local weighted heights. More precisely, given $\x \in \X$, we let $K$ be a finite extension of $k$ containing   $k(\x)$ and define:
\begin{equation}\label{gwh2}
	\hn_{\Y}(\x)\colonequals\sum_{u\in M_K}^{} \il_{ \Y} (\x, u),
\end{equation}
which  is independent of the choice of the field $K$.   

We notice that for any weighted ideal sheaf $\cI \subset k_\w[x_0, \cdots, x_n]$, we may define $\il_\cI$ and 
$ \hn_{\cI}$ as   $\il_\Y$  and  $ \hn_{\Y}$ where $\Y=\Zc(\cI) $ is the weighted variety define by $\cI$.
The basic properties of weighted local and global heights associated to closed subschemes  are respectively described in Propositions 4.6 and 4.7 of \cite{Salami-Shaska1}.

%The   weighted  local and global  heights associated to closed subschemes of  weighted projective varieties are  in  \cite[Subsection 4.5]{Salami-Shaska1}.

\subsection{Weil and Cartier divisors on weighted varieties}
Let $\X$  be a weighted  variety in $\P^n_{\w, k}$ over the field $k.$ 
The group of \textbf{Weil divisors} on $\X$  is a free Abelian group generated by weighted closed  subvarieties of codimension one on $\X.$  This group is denoted by  $\WeDiv_\w   (\X)$. 
The \textbf{support} of a Weil divisor  $D=\sum_Y n_\Y \cdot \Y $  is the union of all  $\Y$~'s  such that $n_\Y \neq 0$, which is denoted by $\Supp   (D)$.
Such a   divisor is said to be \textbf{effective} if every $n_\Y \geq 0$ for all codimension one subvarieties $\Y \subset \X$.    

We define $\ord_\Y: \cO_{\X,\Y} \setminus \{0\} \rightarrow \Z$  to be 
$$\ord_\Y(f)=\text{length}_{\cO_{\X,\Y}} \left(\frac{\cO_{\X,\Y}}{ \left\langle f\right\rangle } \right), $$
which is well defined since $\cO_{\X, \Y}$ is a local ring. 
Then, one can  extend $\ord_\Y$ to the fraction field $k_\w(\X)^*$ in the usual way.
%
%\begin{lem}
The order function $\ord_\Y: k_\w(\X)^* \rightarrow \Z$ has the following properties:
\begin{enumerate}[\upshape(i)]
	\item  $\ord_\Y(f\cdot g)=\ord_\Y (f)+\ord_\Y (g)$
	\item  For a fixed $f \in k_\w (\X)^*$ there are only finitely many $\Y$ such that $\ord_\Y \neq 0$.
	\item Let $f\in k_\w(\X)^*$. Then,  $f \in \cO_{\X,\Y}$ if and only if $\ord_\Y (f)\geq 0$.
	Similarly, $f \in \cO_{\X,\Y}^*$ if and only if $\ord_\Y (f)=0$.
	\item If $\X$ is weighted projective  variety and $f \in k_\w(\X)^*$, then  $f \in k^*$ if and only if  $\ord_\Y (f) \geq 0$ for all $\Y$; if and only if  $\ord_\Y (f)=0$ for all $\Y$.
\end{enumerate}
The Weil divisor of any  rational function  $f \in k_\w(\X)^*$ is defined as 
\[ \div (f)=\sum_{\Y \subset \X } \ord_\Y (f)\cdot \Y \]
which is called a \textbf{principal divisor} on  weighted  variety $\X$.
Two Weil divisors $D$ and $D'$ on a weighted variety $\X$ are said to be \textbf{linearly equivalent} if their difference is a principal divisor.
The  Weil divisor of zeros and  divisor of poles  of $f$,     are  respectively defined as:
\[ (f)_0=\sum_{\ord_\Y(f)>0} \ord_\Y (f) \cdot \Y,\ \text{and} \  \ (f)_\infty=-\sum_{\ord_\Y <0} \ord_\Y (f) \cdot \Y\]
The \textbf{divisor class group} of $\X$ is the group of divisor classes modulo linear equivalence. This group is denoted by $\Cl_\w   (\X)$,  and  $\Cl ( \P^n_{\w,k})$  for $\X =  \P^n_{\w,k}$.

Now, we define a \textbf{Cartier divisor} on a weighted variety $\X$ as an equivalence class of collection of pairs 
$(U_i, f_i)_{i\in I}$   such that
the followings hold:
\begin{enumerate}[\upshape(i)]
	\item The $U_i$ are affine weighted open  sets that cover $\X$.		
	\item The $f_i$ are non zero rational functions, $f_i \in k_\w(U_i)^*=k_\w(\X)^*$.		
	\item $\frac{f_i}{f_j} \in \cO_\X{(U_i\cap U_j)}^*$, so $\frac{f_i}{f_j}$ has no poles or zeros on $U_i\cap U_j$.		
\end{enumerate}

\noindent Two Cartier divisors $\{(U_i, f_i)| i \in I\}$ and $\{(V_j, g_j)| j \in J\}$ on a weighted variety $\X$ are said to be  equivalent if  for all $i \in I$ and $j \in J$ one  has 
\[
\frac{f_i}{g_j}\in \cO_\X(U_i\cap V_j)^*.
\]
The \textbf{sum of two Cartier divisors} is defined in as follows:
\[\{(U_i, f_i)| i \in I\}+\{V_j, g_j)| j \in J\}=\{(U_i\cap V_j , f_i \, g_j)| (i,j) \in I\times J \}.\]  
The set of  Cartier divisors with this operation on a weighted variety $\X$ form a group that we denote it by $\CaDiv_\w (\X)$.  The \textbf{support} of a Cartier divisor is the set of zeros and poles of  $f_i$, which is denoted by $\Supp(D)$.  A Cartier divisor is said to be \textbf{effective} or \textbf{positive} if it can be defined by a collection $\{(U_i, f_i)| i \in I\}$  such that every $f_i \in \cO_\X(U_i).$  
For a given  $f\in  k_\w(\X)^*$, the divisor  $\div(f)=\{(\X,f)\}$ is called a \textbf{principal Cartier divisor}.  Two Cartier divisors  are \textbf{linearly equivalent}  if their difference is a principal divisor.  The group of Cartier divisors classes modulo linear equivalence is called \textbf{Picard group} of a weighted variety $\X$, which is denoted by $\Pic_\w (\X)$, and   $\Pic(\P_{\w ,k}^n)$
when $\X= \P_{\w, k}^n $. 

For reduced vector of weights $\w$, 	 in \cite[Sections 5, 6]{AlAmrani1989},
it is proved  that  the  maps $\Z \rightarrow \Cl(\P_{\w ,k}^n), $ given by $	1 \mapsto \cO_\X (1),$ and 
$	\Z \rightarrow \Pic(\P_{\w ,k}^n),$   given by $ 	1 \mapsto \cO_\X (m), $ with $ 	\ m=\lcm(q_0, \cdots, q_n)$
induce    the isomorphisms $\Cl(\P_{\w ,k}^n ) \iso \Z, $ and $ \Pic( \P_{\w ,k}^n ) \iso \Z, $  respectively.
Moreover,  $\cO_{\P_{\w ,k}^n} (a)$ is not necessarily an invertible sheaf for any given integer $a\in \Z$. 
However,   the sheaf $\cO_{\P_{\w ,k}^n } (m)$  with $ m=\lcm (q_0, \cdots, q_n)$ 
is  ample and invertible, and for $a, b \in \Z$ we have 
\[
\cO_{\P_{\w ,k}^n} (a) \otimes \cO_{\P_{\w ,k}^n} (m)^{\otimes b} \iso \cO_{\P_{\w ,k}^n} (a+b m).
\]
Furthermore,   it is proved that  $ \cO_{\P_{\w, k}^n} (m)$ 
is ample  and there is $c \in \Z$  depending on $q_i$'s   such that  $ \cO_{\P_{\w, k}^n} (c m)$ is very ample.
For more detail on the proofs of above assertions, one can see  by  \cite[Thm. 4B.7 in Sec. 4]{Beltrametti1986}.
For example, when $\w=(q_0, q_1, q_2)$ is reduced and well-formed, then  $\cO_{\P_{\w ,k}^2} (m)$ with $m=q_0 q_1 q_2$	 is both ample and  very ample. For  $\w=(1, 6, 10, 15)$  and $m=30$, the invertible   sheaf   $\cO_{\P_{\w ,k}^3} (30)$ in only ample, but   
$\cO_{\P_{\w ,k}^3} (30 \ell )$  is very ample for all $ \ell \geq 2$.

%	Furthermore, the sheaf  $ \cO_{\P_{\w, k}^n} (a)$ is   coherent  and Cohen-Macaulay for any  $a \in \Z$.
%	If  $ \cO_{\P_{\w, k}^n} (a) \neq 0, $ then it is reflexive of rank $1$ by  \cite[Cor. 5.8]{Beltrametti1986}.

For a  subset $J=\{ i_0, \cdots, i_r\} \subset \{ 0,1, \cdots, n\}$  with $r\leq n$ and 
a vector of weights $\w=(q_0,  \cdots, q_n)$, we set
$\w_J=(q_{i_0},  \cdots, q_{i_r})$  and $m_J= \lcm \{ q_{i_0},  \cdots, q_{i_r}\}$.
A Cartier divisor $D$ on a weighted variety $\X \subseteq \P^n_{\w, k}$ is said to be \textbf{very ample} if  there exist a closed immersion $X 	\hookrightarrow  \P^n_{\w_J, k}$ for some subset $J=\{ i_0, \cdots, i_r\} \subset \{ 0,1, \cdots, n\}$ such that 
$\cO_{\X} (D)  \iso   \cO_{\P_{\w_J ,k}^r} (m_J) |_\X$.
We say that  a  Cartier divisor $D$ on a weighted variety $\X $ is 
\textbf{ample}   if $\ell_J D$ is very ample  for some positive large integer $\ell_J$.
In the 	rest of  paper, we apply the same terminologies to the associated invertible sheaf or line bundle.
% $\cL_D=\cO_{\X} (D)$.

\subsection{Vojta's conjecture for smooth  weighted varieties}
We denote by  $\sing(\P_{\w}^n )$ the singular locus  of $\P_{\w}^n$.
If $\left[x_{0}: \ldots: x_{m}\right] \in \sing(\P_{\w}^n )$, then $x_{i}=0$ for at least one $i$. 
%  \tony{We have denoted points in $\P_{\w}^n$ by $[x_0, \dots , x_n]$ not by tuples $\left(x_{0}: \ldots: x_{m} \right)$.}
Moreover, if  $\gcd \left(a_{j}, a_{j}\right)=1 $ for every couple   $(i, j)$ with  $j \neq i$, then
$ \sing(\P_{\w}^n)=\left\{P_{0}, \ldots, P_{m}\right\}$ 
where $P_{i}$ are the $m+1$ vertices $[0: \ldots: 1: \ldots: 0]$.

%Next we consider the  canonical  quotient map  $p_\w:  \V_k^{n+1}\rightarrow \P_{\w, k}^n $,  which induces the surjective morphism 
%$\pi_\w: \P_k ^n  \rightarrow \P_{\w, k}^n$.
%
%Let $\X$ be a weighted subvariety  of  $\P_{\w, k}^n $.  
A  weighted subvariety   $\X$  of  $\P_{\w, \Q}^n $ is called \textbf{quasi-smooth}  of dimension $m$   if its affine cone  $\cC_\X$ is smooth   variety of dimension $m+1$ outside its vertex.     The singularities of a quasi-smooth variety $\X$  are due to the $k^\ast$-action  and hence are cyclic quotients singularities.
Furthermore, by \cite[Cor. 5.9]{Beltrametti1986}, if $\X \subset \P^n_{\w, k}$ is subvariety such that
$\X \cap \sing (\P^n_{\w, k})= \emptyset$, then  $\X$ is non-singular if and only if $\X$ is quasi-smooth.
%Recall from \cite{Salami-Shaska1} that  $\P_{\w, \Q}^n $ is regular with codimension one and if $ q_i>1$'s  are mutually coprime  then  
%
%\[\sing(\P_{\w, \Q}^n )=\{\x_i=[0: \cdots: 1: \cdots: 0] : \  0\leq i \leq n \}.\]
%
%\begin{defi} \label{wps-smooth}
%A non-singular weighted projective space $\P^n_{\w, \Q} $ is  called 

%\end{defi}

%	In the rest of paper, a {\bf smooth weighted projective variety} $\P^n_{\w, \Q} $ means a   quasi-smooth weighted projective variety with $\X \cap \sing (\P^n_{\w, k})= \emptyset$.
Let $\X$ be a smooth weighted projective variety  in $\P^n_{\w, \Q}$
defined  over   $\Q$.
Suppose   that there is an open  subvariety $\U \subset \X$ with  complement of codimension at least two which is Gorenstein, i.e.,   	the  dualizing  sheaf	$\omega_\U^0$ is invertible. We let $K_\U$ be a canonical divisor on $\U$ such that  $\omega_\U^0 \iso \cO_\U(K_\U)$, and then define the  	\textbf{canonical sheaf } $K_\X$ as  the  closure of $K_\U$ in $\X$.     

By  a \textbf{ normal crossing divisor} on weighted $\X$, we mean a Cartier divisor $D$ such that at every point	in the support of $D$  such that $D$ is given locally by an equation of the form  $x_0^{1/q_0} x_1^{1/q_1} \ldots  x_n^{1/q_n}=0$. 
Then,  we   formulate  the analogue of  Conj. \ref{voj1}, for the smooth weighted projective varieties as follows:

\begin{conj}\label{voj1wsp}
	Let $\X$ be a smooth weighted projective variety, $K_\X$ a canonical divisor, $\Ac$ an ample divisor, and $D$ a normal crossings divisor on $\X$, all defined over $\Q$.  Let $S$ be a finite subset of places of $\Q$   containing $M_\Q^\infty$.  Then, given any real constant  $\varepsilon >0$ and any positive integer $r$, there exists   a proper weighted Zariski-closed subset  $Z$ of $\X$, depending only on $ \X,  D,   \Ac,  \varepsilon,  r$, such that 
	\[\hn_{K_\X}(\x) + \sum_{\nu \in S}^{ } \il_D(\x, \nu) \leq  \varepsilon \cdot \hn_\Ac(\x ) +
	% \frac{1}{m}
	d (\x) + O(1),\]
	for all $\x \in (\X \setminus Z) (\bar \Q)$ with $[\Q (\x):\Q ]\leq r$.
\end{conj}

%^^^^^*
\subsection{Vojta's conjecture using a weighted multiplier ideal sheaf}
In order to avoid   using   the weighted normal crossing divisors and replace it by an error term as in \cite{Vojta2012}, 
let $\cI$ be a nonzero  weighted ideal sheaf on a weighted projective variety $\X$ and $c\geq 0$ a real constant.
For the definition of a  weighted ideal, one can see \cite{Salami-Shaska1}.  
Let $f: \X' \rightarrow \X$  be a proper birational morphism such that $\X'$ is smooth weighted variety and  $f^* \cI = \cO_{\X'} (-E)$,  for a normal crossing  divisor $E$ on $\X'$. Denoting by $\Rc_{\X'/\X}$ the ramification divisor  of $\X'$ over $\X$, we define the \textbf{weighted multiplier ideal sheaf} $\cI_c$ and  $\cI_c^-$ associated  to $\cI$ and $c$ as 
\begin{equation}
	\cI_c\colonequals  f_* \cO_{\X'} (\Rc_{\X'/\X}  - \lfloor c E \rfloor),\ \text{and}\  \cI_c^- \colonequals\lim_{\varepsilon \rightarrow 0^+} \cI_{c-\varepsilon}.
\end{equation}
As in the case of projective varieties, we denote $\cI_1$ and  $\cI_1^-$ by  $\cI$ and  $\cI^-$, respectively.

An analogue  of   Conj. \ref{voj2}, for weighted projective varieties follows: 

\begin{conj}\label{wpv-voj2}
	Let $\X$ be a smooth weighted projective variety, $K_\X$ a canonical divisor, $\Ac$ an ample divisor and $\cI$  a nonzero weighted ideal sheaf $\X$ all defined over $\Q$.  Let $S$ be a finite subset of places of $\Q$ containing $M_\Q^\infty$.  Then, given any real constant  $\varepsilon >0$ and positive integer $r$, there exists a proper Zariski-closed subset  $Z$ of $\X$, depending only on $ \X,  \cI,   \Ac,  \varepsilon,  r$, such that 
	\[
	\lwh _{K_\X}(\x) + \sum_{\nu \in S}^{ } \il_\cI(\x, \nu) -  \sum_{\nu \in S}^{ } \il_{\cI^-}    (\x, \nu) \leq  \varepsilon \lwh_\Ac(\x) + 
	%\frac{1}{m} 
	d(\x) + O(1), \]
	for all $\x \in (\X \setminus Z) (\bar \Q)$ with $[\Q(\x): \Q]\leq r$.
\end{conj}

%^^^^^^^^^^
%
\subsection{Vojta's conjecture for weighted log pairs}
Next we follow closely the terminology of the log pairs for projective varieties as in section \ref{log-pairs}.  
A \textbf{$\Q$-divisor}  on a weighted  projective   variety $\X$ is a formal finite sum 
\[
D=\sum c_i D_i.
\]
%	where $c_i\in \Q$   and $D_i\in  \CaDiv_\w(\X)$, where  $\CaDiv_\w(\X)$ denotes the set of Cartier divisor on a weighted variety $\X$,  see \cite[Sec. 4]{Salami-Shaska1} for more details.

A $\Q$-divisor $D$ on weighted   variety $\X$ is said \textbf{integral} if  all coefficients $c_i$'s are integers.	
By clearing the denominators of $c_i$'s, we can write $D= c D'$ for some   $c\in \Q$ and an integral weighted divisor $D'$. 
It is called \textbf{effective} if $c_i\geq 0$      and $D_i$ are weighted integral divisors on $\X$.  
The  \textbf{support}  of $D$,  denoted by $\Supp (D)$, is 
\[
\Supp (D) = \bigcup \Supp (D_i)
\]
as in the case of projective varieties.  $D=\sum c_i D_i$ is  called \textbf{ample} if  $c_i\in \Q$, $c_i> 0$   and $D_i$ are all ample Cartier divisors on $\X$.   Here, by a \textbf{Weil $\Q$-divisor} on a weighted variety $\X$, we mean a $\Q$-linear combination of its codimension one subvarieties, i.e, an  element of  %
\[
\WeDiv_\w (\X) \otimes \Q.
\]
We use   $\lceil D \rceil$ and $\lfloor D \rfloor$ to denote the round up and round down of any    Weil $\Q$-divisor  $D=\sum_{i}^{} c_i \Y_i  $, that is, 
\[
\lceil D \rceil =\sum_{i}^{} \lceil  c_i\rceil \Y_i \quad  \text{    and     } \quad    \lfloor D \rfloor= \sum_{i}^{}  \lfloor  c_i\rfloor \Y_i.
\]
A weighted projective variety $\X$ defined over $\Q$ is said to be \textbf{$\Q$-Gorenstein} if it  is  Gorenstein in codimension one,   satisfies Serre's condition $S_2$,   and a canonical divisor $K_\X$ is $\Q$-Cartier. For example, if $\X$ is normal, then the first two conditions are true   and hence a weighted canonical divisor exists unequally up to linear equivalence  and  is Cartier in  codimension one.
A \textbf{$\Q$-subscheme} of  $\X$ is a formal linear combinations 
\[
\Y=\sum_{i=1}^m c_i \cdot \Y_i
\]
of proper closed subschemes $\Y_i \subset \X$ with  all $c_i\in \Q$. The support  of such $\Y$ is defined to be the closed subset $\displaystyle \cup_{c_i \neq 0} \Y_i$, and  it is called   \textbf{effective $\Q$-subscheme} if $c_i \geq 0$ for every $i$.

By a \textbf{resolution} of a weighted variety  $X$ over $\Q$, we mean  a projective birational morphism
$f: \X' \rightarrow \X$ such that   $\X'$ is a smooth projective weighted variety  over $\Q$.   
The existence of such  a weighted resolution is proved by Dolgachev in \cite[Subsec. 4.2]{dolga}, by giving a construction of resolution for the sheaf of regular differentials on $\X$.  	One can see \cite{Wlod2023}, for a Rees algebra approach on weighted resolution.

A \textbf{weighted log pair} is called  a pair $(\X, \Y)$ of a  $\Q$-Gorenstein weighted variety $\X$  and an effective  weighted $\Q$-subscheme $\Y$ of $\X$. 	A \textbf{weighted log resolution} of  a weighted log pair $(\X, \Y)$  with 
\[
\Y=\sum_{i=1}^m c_i \cdot \Y_i,
\]
is a  	projective birational morphism  $f: \X' \rightarrow \X$ of $\X$ such that
$\X'$ is a weighted smooth variety defined   over $\Q$,  
the set-theoretic inverse image $f^{-1} (\Y_i)$  is a weighted Cartier divisor on $\X'$,
and the union of $\Exc(f) $ of the exceptional divisor  of $f$ with all $f^{-1} (\Y_i)_{\red}$   is a simple weighted normal crossing divisor of $\X'$.   
The existence of a  the weighted log resolution  of   $(\X, \Y)$ is a consequence of  Hironoka's theorem \cite[Thm. 4.1.3]{Lazarsfeld2004c}.

For a weighted log resolution  $f: \X' \rightarrow \X$  as a weighted log pair $(\X, \Y)$, the \textbf{relative weighted canonical divisor}  of $\X'$ over  $(\X, \Y)$ is defined to be the  weighted $\Q$-Weil divisor 
\[
K_{\X'/(\X, \Y) }=K_{\X'/\X } - f^* \Y, 
\]
where  $f^* \Y$ is the pull-back of $\Y$ by $f$ over $\X'$ and $K_{\X'/\X } $ is the  relative canonical divisor  of $\X'$ over $\X$, i.e., 
\[
\cO(K_{\X'}) \iso \cO(K_{\X'/ \X}) \otimes f^* \cO_\X(K_\X),
\]
which is a $\Q$-Weil divisor on $\X'$.
Given a  weighted log pair $(\X, \Y)$ and  a  weighted log resolution  $f: \X' \rightarrow \X$, we write 
\[
K_{\X'/(\X, \Y) } =\sum_{\Z}^{} a_\Zc \cdot \Zc,
\]
where $\Zc$ runs over all prime  divisors of $\X^\prime$. The weighted log pair $(\X, \Y)$ is called
\textbf{strongly canonical } (resp.
\textbf{Kawamata log terminal,} and  \textbf{log canonical}) if $a_\Zc \geq 0$ (resp. $a_\Zc >0$, and $a_\Zc \geq -1$ ) for every $\Zc.$
These properties are independent of the resolution and are also local.

Define the \textbf{weighted non-strongly canonical locus (wnsc)}   of  the pair $(\X, \Y)$ to be the  smallest weighted closed subset  $\Wc \subset \X$ such that the  weighted pair $(\X \backslash \Wc , \Y|_{\X \backslash \Wc})$   is    strongly canonical. 
The  \textbf{weighted non-Kawamata log terminal locus (wnklt)}   of  the pair $(\X, \Y)$ to be the  smallest weighted closed subset  $\Wc \subset \X$ such that the  weighted pair is  Kawamata log terminal.
Similarly the  \textbf{weighted non-log canonical (wnlc) })  of  the pair $(\X, \Y)$ to be the  smallest weighted closed subset  $\Wc \subset \X$ such that the  weighted pair is       log canonical.  
We denote them respectively  as $\wnsc (\X, \Y)$,        $\wnklt (\X, \Y)$,  $\wnlc (\X, \Y)$. One may check that 
\[
\wnsc (\X, \Y) \subset  \wnklt (\X, \Y)  \subset \wnklt (\X, \Y).
\]
For a weighted log pair $(\X, \Y)$ with a log resolution $f: \X' \rightarrow \X$, we define $\cI(\X, \Y)$ a \textbf{variant of multiplier sheaf} as
\[
\cI(\X, \Y)\colonequals f_* \cO_{ \X'} (\lceil K_{\X'/(\X, \Y) } \rceil)
\]
if $\X$ is a weighted normal variety; see \cite[9.3.56]{Lazarsfeld2004b} for the definition of multiplier ideal sheaf in usual case.      Otherwise, we set
\[
\cI(\X, \Y)\colonequals \bar{f}_* \cO_{ \X'} (\lceil K_{\X'/(\X, \Y) } \rceil),
\]
where $\bar{f}_*   \cO_{ \X'}(E)$ denotes the largest  ideal sheaf  in $\cO_{ \X}$ for which its pull-back by $f$ is contained in  $\cO_{ \X'}(E)$ as an $\cO_{\X'}$-submodule of   (constant) function field sheaf $\M_{\X'}$. 
Moreover, there exist a constant $\varepsilon_0>0$ such that for every rational number    $0< \varepsilon \leq \varepsilon_0$, one has  
$\cI(\X, (1-\varepsilon) \Y)=\cI(\X, (1-\varepsilon_0) \Y)$.

Let  
\begin{equation}
	\begin{split}
		& \cI^-(\X, \Y)\colonequals\cI(\X, (1-\varepsilon) \Y), \  (0< \varepsilon \ll 1) \\
		&   \cH(\X, \Y)\colonequals\bar{f}_* \cO_{ \X'} (\lfloor K_{\X'/(\X, \Y) } \rfloor), 
	\end{split}
\end{equation}
where $\bar{f}_* $  is as above.    We note that the  definition of $\cH(\X, \Y)$ does not depend on the notion of  ``weighted simple normal crossing".

\begin{lem}
	The definitions of $\cI(\X, \Y)$,   $\cI^-(\X, \Y)$, and $\cH(\X, \Y)$ are  independent of the choice of a weighted log resolution.
\end{lem}

\proof
By  adopting \cite[Lem. 3.1]{Yasuda2018} and  \cite[Prop. 3.4]{Yasuda2018} respectively to the case of weighted projective schemes,
one get the result for  $\cI(\X, \Y)$ and    $\cI^-(\X, \Y)$. An argument similar to the proof of the Proposition 3.4 in   \cite{Yasuda2018} shows the  assertion for $\cH(\X, \Y)$. 
\qed

\begin{prop}\label{Wlocus}
	Let   $(\X, \Y)$ be a weighted log pair. Then, the following are true: 
	\begin{enumerate}[\upshape(i)]
		\item  $\wnlc (\X, \Y)\subset \Supp (\cO_{ \X}/\cI^-(\X, \Y) ) \subset   \wnklt (\X, \Y);$
		
		\item If $(\X \backslash \Supp (\Y), 0)$ is a weighted Kawamata log terminal, then
		\[\wnlc (\X, \Y)= \Supp (\cO_{ \X}/\cI^-(\X, \Y) );\]
		
		\item $ \Supp (\cO_{ \X}/\cH(\X, \Y) )= \wnsc (\X, \Y).$

		\item If $(\X, \Y)$  is weighted log canonical, then  $\cO_{ \X}/\cH(\X, \Y)$ is reduced, i.e., as  a reduced closed subscheme,  $\cH(\X, \Y) $   is the   ideal sheaf of $\wnsc (\X, \Y)$.
	\end{enumerate}	
\end{prop}

\proof  Let $f:\X' \rightarrow \X$ be a weighted log resolution of $(\X, \Y)$, and denote by $\mult_\Zc(E)$ the multiplicity of any divisor $E$ on $\X'$.

(i) Given any prime divisor $\Zc$ of $\X'$ and real constant $0< \varepsilon \ll 1$, we have 
\[
\mult_\Zc ( K_{\X'/(\X, \Y)} ) < -1,
\]
which implies that
\[ 
\mult_\Zc ( K_{\X'/(\X, \Y)}+   \varepsilon f^* \Y) < 0, 
\]
and hence
$\mult_\Zc ( K_{\X'/(\X, \Y)} )  \leq -1.$
This proves part (i).

(ii)  It is enough to show  that if  the pair  $(\X \backslash \Supp(\Y), 0)$    is a weighted Kawamata log terminal, then   
\[
\mult_\Zc ( K_{\X'/(\X, \Y)} )  \geq -1
\]
and so 
\[
\mult_\Zc ( K_{\X'/(\X, \Y)}+   \varepsilon f^* \Y) \geq  0. 
\]
If  $\mult_\Zc ( K_{\X'/(\X, \Y)} )>-1$, then the result  is trivial. 
If   $\mult_\Zc ( K_{\X'/(\X, \Y)} )=-1$,   then  $\Zc$ is contained in  $\Supp(f^* \Y)$ by assumption on   $(\X \backslash \Supp(\Y), 0).$ Thus,   
\[
\mult_\Zc ( K_{\X'/(\X, \Y)}+   \varepsilon f^* \Y) >-1  \Rightarrow \mult_\Zc ( K_{\X'/(\X, \Y)}+   \varepsilon f^* \Y) \geq 0.
\]

(iii) If $(\X, \Y)$ is a weighted strongly canonical, then $ \lfloor K_{\X'/(\X, \Y)} \rfloor \geq 0.$ 
By definition of   $\bar{f}$, we have  
\[
\bar{f}_* \cO_{ \X'} (\lfloor K_{\X'/(\X, \Y) } \rfloor)=\cO_\X.
\]
This shows that 
\[
\wnsc (\X, \Y)  \subset \Supp (\cO_{ \X}/\cH(\X, \Y) ).
\]
If $(\X, \Y)$ is not a weighted strongly canonical around $x\in \X$, then there is a prime divisor $\Zc$ on $\X'$
such that $x\in f(\Zc)$ and $\mult_\Zc ( K_{\X'/(\X, \Y)} ) <0. $ Thus, we have 
\[
\cO_{\X'} \not \subset \cO_{\X'} (\lfloor K_{\X'/(\X, \Y) } \rfloor).
\]
Replacing $\X$ with any open neighborhood of $x$ does not change  the last result. Therefore, 
\[
\Supp (\cO_{ \X}/\cH(\X, \Y) )  \subset   \wnsc (\X, \Y).
\]

(iv) By part (iii), we have $\cH(\X, \Y)$ is a subset of the ideal sheaf of $\wnsc (\X, \Y)$ denoted by $\Nc$.
Now, let $\U \subset \X$  be an open set and $g \in \Nc (\U)$. Then $f^* g$  vanishes  along the closed set 
$f^* ( \wnsc (\X, \Y))$   containing every prime divisor $\Zc$ on $\X'$ having negative coefficient in
$\lfloor K_{\X'/(\X, \Y) } \rfloor$, which is equal to $-1$ since $(\X, \Y) $ is weighted log  canonical.
Therefore, 
\[
f^* g \in \cO_{\X'} (\lfloor K_{\X'/(\X, \Y) } \rfloor) (f^{-1} \U),
\]
which implies $g\in \cH (\U)$ and hence $\cH \subset \Nc.$
\qed

Given a Weil $\Q$-divisor  $D$ on a weighted projective variety $\X$,  such that $n D$ is a Weil divisor on $\X$, we define
weighted  height function 
\[
\lwh_D:= \frac{1}{n} \lwh_{nD}.
\]

For a weighted $\Q$-Gorenstein projective variety  $\X$ with a canonical divisor $K_\X$, we can define  a weighted global height  function $\lwh_{K_\X}$ up to addition of a bounded function. 
Given a weighted  log pair $(\X, \Y)$ of a weighted $\Q$-Gorenstein $\X$,   we define  the weighted height function  associated to the subscheme $K_{(\X, \Y)}=K_{\X}+ \Y$ as  
\begin{equation}\label{eqsum}
	\lwh_{K_{(\X, \Y)}}=\lwh_{K_\X}+ \lwh_\Y,
\end{equation}
where $\lwh_\Y$ is the weighted height function associated to the subscheme $\Y$ or its  ideal sheaf.
For more detail on the local and global weighted height functions, we refer the reader to \cite{Salami-Shaska1}.

Next we are ready to state  Vojta's conjecture for weighted log pairs.

\begin{conj}\label{voj6}
	Let   $\X$ be a weighted projective scheme, $\Y$ a closed weighted subscheme with ideal sheaf $\cI=\cI(\Y)$, 
	$K_\X$ a canonical divisor, 
	and  $\Ac$ an ample divisor on $\X$,  all defined over  $\Q$. 
	
	Let  $(\X, \Y)$ be a weighted log pair and  $S$ be a finite subset of places of $\Q$ containing $M_\Q^\infty$.  Then, given any real constant  $\varepsilon >0$ and a positive integer $r$, there exists  a proper weighted  Zariski-closed subset  $Z$ of $\X$, depending only on $  \X,  \cI,   \Ac,  \varepsilon,  r$, such that 
	\[
	\lwh_{K_{(\X, \Y)}  (\x)} - \sum_{\nu \not \in S}^{ } \il_\cH(\x, \nu) -  \sum_{\nu \in S}^{ } \il_{\cI^-}(\x, \nu) \leq  \varepsilon \lwh_\Ac(\x) +  
	%\frac{1}{m}
	d  (\x) + O(1).
	\]
	for all $\x \in (\X \setminus Z) (\bar \Q)$ with $[\Q (\x): \Q]\leq r$, where 
	\[
	\cH=\cH(\X, \Y) \quad \text{ and } \quad \cI^-=\cI^-(\X, \Y).\]
\end{conj}

We note that the terms 
\[
\sum_{\nu \in S}^{ } \il_\cH(\x, \nu) \quad \text{  and  } \quad \sum_{\nu \in S}^{ } \il_{\cI^-}(\x, \nu)
\]
can be thought of as the contribution of   $\wnsc (\X, \Y) $ and $\wnklt (\X, \Y)$,   or $\wnlc (\X, \Y) $ if $(\X \backslash \text(Supp) (\Y) , 0)$  is Kawamata log terminal.

Since a pair $( \X, D)$ with $\X$ a smooth weighted variety and $D$ a reduced simple weighted normal crossing divisor on $\X$ is a weighted log  canonical, by parts (i), (iii) and (iv) of   Prop. \ref{Wlocus}, one can conclude  that  
\[
\sum_{\nu \in S}^{ } \il_{\cI^-}(\x, \nu) =0
\]
and hence  the right hand side of  the inequality of  Conj. \ref{voj6} is equal to  
\[
\lwh_{K_{(\X, D)}(\x) } - \sum_{\nu \not \in S}^{ } \il_\cH(\x, \nu) =\hn_{K_\X}(\x) + \sum_{\nu \in S}^{ } \il_D(\x, \nu).
\]
Thus  Conj. \ref{voj6}  is the same as   Conj. \ref{wpv-voj2} and   Conj. \ref{voj1wsp},  in this case.

In contrast,  given  a weighted   log pair  $(\X, \Y)$  and a log resolution
\[
f: \X'\rightarrow \X,
\]
if we suppose that   Conj. \ref{voj1wsp}, holds  for $\X'$ and  the reduced simple normal crossing divisor
\[ \lceil K_{\X'/(\X, \Y)}  +\varepsilon f^* \Y \rceil  - \lfloor  K_{\X'/(\X, \Y)} \rfloor,\]
for $0 < \varepsilon \ll 1$, then    Conj. \ref{wpv-voj2} and   Conj. \ref{voj6}   holds for  $(\X, \Y)$.

Indeed, the argument is similar to those given in   \cite[Prop. 4.3]{Vojta2012} and \cite[Prop. 5.4]{Yasuda2018} as follows.
By definition, we have 
\[
f^{-1} \cH  \subset \cO_{\X'} (\lfloor  K_{\X'/(\X, \Y)} \rfloor)   
\]
and 
\[   
f^{-1} \cI^{-} \subset  \cO_{\X'}  ( \lceil K_{\X'/(\X, \Y)}    +\varepsilon f^* \Y \rceil )
\]
for  $0 < \varepsilon \ll 1$. 
Using these and the properties of weighted height functions, we get
\[
\begin{split}
	& \il_\cH \circ f  \geq \il_{-D_1},  \ \text{with } \ D_1= \lfloor  K_{\X'/(\X, \Y)} \rfloor,    \\
	& \il_{\cI^-} \circ f  \geq \il_{-D_2} \ \text{with } \ D_2= \lceil K_{\X'/(\X, \Y)}   +\varepsilon f^* \Y \rceil.
\end{split}
\]
Then, using
% the fact $ \lwh_{K_{\X'}} = \lwh_{K_{(\X, \Y)}}  + \lwh_{K_{\X'/(\X, \Y)}}$   and
the above inequalities, we have
\begin{align*}
	\left( \lwh_{K_{(\X, \Y)}  } - \sum_{\nu \not \in S}^{ } \il_\cH(\cdot, \nu) -  \sum_{\nu \in S}^{ } \il_{\cI^-}(\cdot, \nu)\right) \circ f &   \leq \lwh_{K_{\X'}}  -  \lwh_{K_{\X'/(\X, \Y)}} 
	- \sum_{\nu \not \in S}^{ } \il_{-D_1}  -  
	\sum_{\nu \in S}^{ } \il_{-D_2 } \\
	&  \leq \lwh_{K_{\X'}}  +  \lwh_{-D_1 } - \sum_{\nu \not \in S}^{ } \il_{-D_1}  -  
	\sum_{\nu \in S}^{ } \il_{-D_2 }	\\
	&  \leq \lwh_{K_{\X'}}  +  \sum_{\nu \not \in S}^{ } \il_{D_2-D_1},
\end{align*}
where 
\[
D_2-D_1=\lceil K_{\X'/(\X, \Y)}   +\varepsilon f^* \Y \rceil -\lfloor  K_{\X'/(\X, \Y)} \rfloor.
\]
Recall that the pullback $f^*  D$ of an ample divisor $D$ is ample.  Therefore,
the above argument leads to  following result.

\begin{cor}\label{eq-conj}
	The Conjecture \ref{voj1wsp},    Conjecture  \ref{wpv-voj2},  and Conjecture  \ref{voj6} are equivalent.
\end{cor}

%^^^^^^^^^^^^^^^^^^^^^^^^^^^
\section{Weighted blow-ups and generalized weighted gcds}\label{sect-6}
%
% In \cite{MR2162351} Silverman used the idea of generalized \gcds to define a height for blow-ups of smooth projective varieties and then assuming  Vojta's conjecture for such height function obtained some conjectural   results  on the generalized gcds.    
%
In this section, we  define the generalized weighted greatest common divisor generalizing the usual one. Then, we describe it as
a weighted height  for weighted blow-ups using the theory introduced in \cite{Salami-Shaska1}. 
%As above,  $k$ is a number field, $\cO_k$ its ring of integers, and $\nu_p$ the valuation at a prime $p \in \cO_k$. 

For any two integers $\a, \b \in \Z$ the \textbf{greatest common divisor } is defined as 
\[ 
\gcd (\a, \b) \colonequals \prod_{p \ \text{prime}} p^{\min \{ \nu_p (\a), \, \nu_p (\b) \}}, 
\]
where $\nu_p$ is the place associated to the prime number $p$.
%\dc{If $\cO_k$ is a PID, then we define $\gcd (\a, \b)$ as the generators of the above ideal. }
The \textbf{logarithmic greatest common divisor }  is % $ \log \gcd (\a, \b)$.
\[ \log \gcd (\a, \b) \colonequals  \sum_{p  \text{ prime}}  \min \,  \{ \nu_p (\a), \nu_p (\b) \}  \log p.
%	= \sum_{p\  \text{prime}} \min \,  \{ \nu_p (\a), \nu (\b) \}    
\]
For each place $\nu \in M_\Q$, we define  
\begin{equation}\label{nu+}
	\begin{split}
		\nu^+ : \Q      & \longrightarrow [0, \infty], \\
		\alpha          &  \mapsto  \max \{ \nu(\alpha), 0\},  
	\end{split}
\end{equation}
which  $\nu^+$ can be viewed as a height function on $\P^1  (\Q) = \Q \cup \{ \infty \}$ with respect the divisor $(0)$,  where we set $\nu^+ (\infty) = 0$.
The \textbf{generalized logarithmic greatest common divisor} of  $\a, \b \in \Q$ is defined  as
\[ 
\hgcd (\a, \b) \colonequals \sum_{\nu \in M_\Q} \min \{ \nu^{+} (\a), \nu^{+} (\b) \}. 
\]
%
%; see \cref{nu+}. 
Then,  given $(\a, \b)\neq (0, 0)$, one may consider the following function 
\begin{equation}
	\begin{split}
		G_{\nu_p} : \P^1  (\Q) \times \P^1  (\Q)    &  \to [0, \infty],  \\
		(\a, \b)  &\mapsto \min \{ \nu^+ (\a), \nu^+ (\b) \},
	\end{split}  
\end{equation}
as a \textbf{local height function} and  the   generalized logarithmic greatest common divisor,   being their sum all together, 
\begin{equation}
	\hgcd (\a, \b) = \sum_{\nu_p \in M_\Q} G_{\nu_p}. 
\end{equation}
as a \textbf{global height function} on   $ \P^1  (\Q) \times \P^1  (\Q).$
% 
%In \cite{MR2162351}, it was given a theoretical interpretation of the  function $G_\nu$ in terms of blow-ups. 

Here is a theoretical interpretation of the  function $G_\nu$ in terms of blow-ups.  More precisely, for $\X= (\P_\Q^1)^2$ let  
$\pi: \tX \rightarrow \X$ 
be the blow-up of  the point $(0,0)$ and $E=\pi^{-1}(0,0)$  be the exceptional divisor of the blow-up. Then, for all $(\a, \b) \in \X  (\Q) \backslash {(0,0)}$ and $\nu\in M_\Q$, one has
\[
\lambda_{\X, E}(\pi^{-1}(\a, \b), \nu)= \min \, \{ \nu^+ (\a),  \nu^+ (\b)\},
\]
and adding all together for  $\nu \in M_\Q$ leads to  
\[
\hgcd (\a, \b)= h_{\tX, E}(\pi^{-1}(\a, \b)).
\]
%By this realization, in  \cite[Def. 2]{MR2162351},  Silverman's  introduced    generalized logarithmic greatest common divisors of  points on  smooth varieties with respect to its subvarieties. 

%^^^^^^^^^^^^^^^^*
\subsection{Generalized  weighted greatest common divisors as heights for blow-ups}
The weighted greatest common divisor for any tuple of integers  $(x_0,  \cdots, x_n) \in  \Z^{n+1}$ was defined in \cite{b-g-sh}, which we are going to recall in below.

Let $\tx = (x_0, \dots, x_n ) \in  \Z^{n+1}$  with $r=\gcd(x_0, \dots, x_n)$ and   
\[
r=u \cdot \prod_{j=1}^{s} p_j,
\]
where $u\in \{ -1, 1\}$  and $p_1, \cdots, p_s$ are prime numbers.    
The \textbf{weighted greatest common divisor} of $\tx   \in  \Z^{n+1}$ is defined as
\begin{equation}
	\wgcd (\tx) \colonequals \prod_{\stackrel{p \in \{p_1, \cdots, p_s\}}{p^{q_i} | x_i} }^{} p
	= \prod_{ \nu_p \in  M_\Q}    p^{ \min \left\{ \left\lfloor  \frac {\nu_p (x_0)}  {q_0}   \right\rfloor,  \dots ,  \left\lfloor  \frac {\nu_p (x_n)}  {q_n}   \right\rfloor  \right\}},
\end{equation}
where the last equality comes from \cite[Lem.~4]{b-g-sh}.   
Here,    the symbol  $\left\lfloor \cdot   \right\rfloor$ denotes the integer part function.    A tuple $\tx \in   \Q^{n+1}$ is said to be \textbf{normalized}  if $\wgcd (\tx)=1$. 
In \cite[Lem.~7 and Cor. 1]{b-g-sh}, it is proved that any point $\x $ in a well-formed space $\P^n_{\w, \Q}$ has a unique normalization   $\y=\frac{1}{\wgcd (\tx)} \star \x$.   
For the rest of  paper, a \textbf{normalized point} $\x \in \P^n_{\w, \Q}$ means a point $\x=[x_0: \ldots : x_n]$ with integer coordinates $x_i\in \Z$ such that $\wgcd (x_0, \ldots , x_n)=1$.

Given   $\w=(q_0, \cdots, q_n)$,    we    let    $\w_i=(1, q_i)$ for each $i=0, 1, \cdots n$. The canonical inclusion 
\[
\Q_{\w_i}(x_i )\hookrightarrow \Q_\w (x_0, \cdots, x_n)
\]
implies the rational map $\P_{\w, \Q}^n \rightarrow  \P_{\w_i, \Q}^1 $ given by   
\[ \x=[x_0: \cdots : x_n] \mapsto [1: x_i],\]
which is defined precisely in the complement of $\{ x_i =0\}$ in $\P_{\w, \Q}^n$.  Considering all of these maps, we  have the rational map 
\begin{equation}\label{phi1}
	\begin{split}
		\varphi_{n, \w}: \P_{\w, \Q}^n   &  \longrightarrow \prod_{i= 0}^n \P_{\w_i, \Q}^1,  \\
		\x=[x_0: \cdots: x_n]  &\mapsto  \varphi_{n, \w} (\x):=( [1 : x_0],   [1 : x_1], \cdots,  [1 : x_n]),
	\end{split}  
\end{equation}
which is defined on the open set $\displaystyle  \P_{\w, \Q}^n  \backslash \cup_{i=0}^n  \{ x_i=0\}$.
For each $p \in  \Z$, we define the function $F_{\nu_p}$  as follows,
\[
\begin{split}
	F_{\nu_p}:\prod_{i=1}^n \P_{\w_i, \Q}^1  &  \longrightarrow \N,  \\
	( [1:x_0],  [1:x_1], \cdots,  [1:x_n])  &\mapsto 
	p^{ \min \left\{ \left\lfloor  \frac {\nu_p^+ (x_0)}  {q_0}   \right\rfloor,  \dots ,  \left\lfloor  \frac {\nu_p^+ (x_n)}  {q_n}   \right\rfloor  \right\}}.
\end{split}
\]

The \textbf{generalized weighted greatest common divisor}  of a given tuple   ${\bar x} =(x_0, \dots, x_n)  \in   \Q^{n+1}$ is defined  as
\begin{equation}\label{log-wgcd}
	\displaystyle
	\hwgcd ({\bar x}) \colonequals \prod_{\nu_p \in M_\Q}   p^{ \min \left\{ \left\lfloor  \frac {\nu_p^+ (x_0)}  {q_0}   \right\rfloor,  \dots ,  \left\lfloor  \frac {\nu_p^+ (x_n)}  {q_n}   \right\rfloor  \right\} },
\end{equation}

We define the  \textbf{logarithmic weighted greatest common divisor}  of any tuple of integers   $\tilde{x}=(x_0,  \cdots, x_n) \in  \Z^{n+1}$ as the sum
\begin{equation}
	\log \wgcd  (\tx) \colonequals  \sum_{\nu \in M_\Q^0}
	%\sum_{\stackrel{\nu \in M_\Q^0}{\nu|_\Q =\nu_p}}  
	\min \left\{ \left\lfloor  \frac {\nu (x_0)}  {q_0}   \right\rfloor,  \dots ,  \left\lfloor  \frac {\nu (x_n)}  {q_n}   \right\rfloor  \right\},
\end{equation}
and the  \textbf{generalized  logarithmic weighted greatest common divisor}  of  any tuple  ${\bar x}=(x_0,  \cdots, x_n)\in \Q^{n+1}$ is  defined to be
\begin{equation}
	\log \hwgcd   ({\bar x}) \colonequals \sum_{\nu \in M_\Q} 
	\min \left\{ \left\lfloor  \frac {\nu^+ (x_0)}  {q_0}   \right\rfloor,  \dots ,  \left\lfloor  \frac {\nu^+ (x_n)}  {q_n}   \right\rfloor  \right\}.
\end{equation}
%
%\begin{defi}
% The \textbf{logarithmic absolute weighted greatest common divisor} $\x \in  \cO_\Q^{n+1}$ is given by 
%
%\[
%\hawgcd (\x) = \frac 1 q \,  \sum_{\nu \in M_\Q^0} \min \, \left\{ \left\lfloor \frac {\nu_p (x_0)} {\bar q_0}   \right\rfloor,    \dots ,  \left\lfloor \frac {\nu_p (x_n)} {\bar q_n}   \right\rfloor  \right\}.
%\]
%
%
%The \textbf{generalized logarithmic absolute weighted greatest common divisor} of $\x$ is given  by
%
%\[
%\hawgcd (\x) \colonequals \frac 1 q \,  \sum_{\nu \in M_\Q^0}  \min \, \left\{ \left\lfloor \frac {\nu^+_p (x_0)} {\bar q_0}   \right\rfloor,    \dots ,  \left\lfloor \frac {\nu^+_p (x_n)} {\bar q_n}   \right\rfloor  \right\}.
%\]
%\end{defi}
%
Let us consider  the following positive real-valued function on $\P_\w^n (\Q)$,
\begin{equation}
	\begin{split}
		T_\nu: \prod_{i=0}^n \P_{\w_i, \Q}^1  & \to [0, \infty] \\
		\left( [1 :x_0], [1: x_n], \dots,  [1 : x_n]  \right) & \to \min \left\{ \left\lfloor  \frac {\nu^+ (x_0)}  {q_0}   \right\rfloor,  \dots ,  \left\lfloor  \frac {\nu^+ (x_n)}  {q_n}   \right\rfloor  \right\}. 
	\end{split}
\end{equation}
For any rational point  $\x=[x_0:  \cdots: x_n] \in \P_\w^n (\Q)$, define its \textbf{generalized  logarithmic weighted greatest common divisor} as
\begin{equation}\label{hwgcd}
	\log \hwgcd (\x) = \sum_{\nu \in M_\Q} T_\nu (\varphi_{n, \w}(\x))=
	\sum_{\nu_{p}  \in M_\Q} \min \left\{ \left\lfloor  \frac {\nu_p^+ (x_0)}  {q_0}   \right\rfloor,  \dots ,  \left\lfloor  \frac {\nu_p^+ (x_n)}  {q_n}   \right\rfloor  \right\},
\end{equation}
where $\varphi_{n, \w} $ is defined by \eqref{phi1}.
%\begin{equation}
%	\log \hwgcd  (\x)=\log \hwgcd  (x_0,  \cdots, x_n).
%\end{equation}
%
\iffalse

\tony{
	\begin{lem}
		All of the above are well defined ..........
	\end{lem}
	
	\proof
	
	\qed
}

\fi

Notice that   all points $\x \in  \P_{\w}^n (\Q)$ with $\log \hwgcd (\x)=0$  belong to the  singular locus $\sing(\P_{\w, \Q}^n )$ as shown next.

%^^^^^^^^^^^^^^^^^^^^^^^^*
\begin{prop}\label{sing1}
	Let $\P_{\w, \Q}^n$ be a   well-formed weighted projective space with $\w=(q_0, \cdots, q_n)$ and   
	$\x    \in \P_{\w }^n  (\Q) $.  If  $\log \hwgcd (\x) =0$  then $\x  \in \sing(\P_{\w, \Q}^n )$.
\end{prop}

\proof    Let $m=\lcm (q_0, \cdots, q_n)$,  and define $J(\x)=\{ j: x_j(\x) \neq 0\}$ for any point $\x=[x_0: \cdots: x_n]  \in  \P_{\w}^n (\Q)  $. Given any prime   divisor $p\mid m$, we  define 
\[
S_\w  (p) = \left\lbrace \x \in  \P_{\w, \Q}^n  : \    p \mid q_i \ \text{for all} \ i \in J(\x)  \right\rbrace.
\]
Then  $\sing(\P_{\w, \Q}^n ) = \bigcup_{p \mid m}^{} S_\w (p)$ and for any prime $p$ we have 

% \[
%\sing(\P_{\w, \Q}^n ) = \left\lbrace  \x \in  \P_{\w, \Q}^n  : \ \gcd( q_i:  i \in \{ j: x_j(\x) \neq 0\})  > 1 \right\rbrace,
%\]
%so we have:
%
\begin{align*}
	\x=[x_0: \cdots: x_n] \in S_\w(p) &   \Rightarrow  p\mid q_i, \ \text{for all } i \in J(\x) \\
	&   \Rightarrow \nu_p^+ (x_i)< q_i, \  \text{for all } i \in J(\x) \\
	&  \Rightarrow   \left\lfloor  \frac {\nu_p^+ (x_i)}{q_i}   \right\rfloor  =0,   \  \text{for all } i \in J(\x) \\
	& \Rightarrow  \min \left\{ \left\lfloor  \frac {\nu_p^+ (x_0)}  {q_0}   \right\rfloor,  \dots ,  \left\lfloor
	\frac {\nu_p^+ (x_n)}  {q_n}   \right\rfloor \right\}  =0. 
\end{align*}
If we assume $\log \hwgcd (\x) =0$, then  
\begin{equation}\label{es1}
	\sum_{\nu_p  \in M_\Q} \min \left\{ \left\lfloor  \frac {\nu_p^+ (x_0)}  {q_0}   \right\rfloor,  \dots ,  \left\lfloor  \frac {\nu_p^+ (x_n)}  {q_n}   \right\rfloor  \right\} = 0.
\end{equation}
Thus, for all $\nu_p \in M_\Q$ with  $p \in \Z$,    we have 
\begin{equation}\label{es2}
	\min \left\{ \left\lfloor  \frac {\nu_p^+ (x_0)}  {q_0}   \right\rfloor,  \dots ,  \left\lfloor  \frac {\nu_p^+ (x_n)}  {q_n}   \right\rfloor  \right\}= 0.
\end{equation}
This  implies that $\x  \in S_\w (p) $ for  any  prime $p\mid m$    and hence $\x \in \sing(\P_{\w, \Q}^n )$.
\qed

The following example shows that the Weil logarithmic height function is strictly positive on $\P_\Q^n$.
\begin{exa} Consider the weights $\w=(1, \dots ,1)$.    Then $\P_{\w, \Q}^n = \P_\Q^n$ is the projective space and the weighted height $\wh_\Q$   is simply the projective height $H_\Q$.      Since $m=\lcm (q_0, \ldots , q_n)=1$ then there are no primes dividing $m$ and $\sing \P_\Q^n = \emptyset$.  
	On the other side from \eqref{hwgcd} we have 
	\[
	\begin{split}
		\log \hwgcd (\x) & =  \sum_{\nu_p  \in M_\Q} \min \left\{ \left\lfloor  \frac {\nu_p^+ (x_0)}  {q_0}   \right\rfloor,  \dots ,  \left\lfloor  \frac {\nu_p^+ (x_n)}  {q_n}   \right\rfloor  \right\}    \\
		& = \sum_{\nu_p  \in M_\Q} \min \left\{         \nu_p^+ (x_0)    ,  \dots ,   \nu_p^+ (x_n)   \right\}  \\
		&  \geq  \min \{ \left\{         \nu_\infty^+ (x_0)    ,  \dots ,   \nu_\infty^+ (x_n)   \right\} >0,
	\end{split}
	\]
	since at least one of the coordinates $x_i \neq 0$, 
\end{exa}

\begin{lem} \label{blow-div}
	Let $\X $ be a smooth weighted variety   defined over $\Q$, and $\Y $ a subvariety of $\X $ of codimension $r\geq 2$.  
	Let  $\pi : \tX \to \X$ be the  blow-up of $\X$ along $\Y$ and denote by  $\tY\colonequals\pi^{-1} (\Y)$  its the exceptional divisor. Then, 
	\begin{enumerate}[\upshape(i)]
		\item  $\pi \, |_{\pi^{-1} (\X\setminus \Y)} \, :  \;  \pi^{-1} (\X\setminus \Y)   \to \X \setminus \Y $ is an isomorphism. 
		\item Exceptional divisor $ \tY$    is an  effective Cartier divisor on $\tX$.
	\end{enumerate}
\end{lem}

\proof   This is a direct consequence of     \cite[Prop. II.7.13]{Ha}.
%
%Since blowing up commutes with restrictions to open subvarieties, the first part means that $\tX = \X$ if $\Y = \emptyset$. Then we are blowing up the ideal $\cI=\cO_\X$ and the result follows. 
% 
%The proof that $\pi^{-1} (\Y)$  is an  effective Cartier divisor on $\tX$ goes similarly to the proof for projective varieties. 
For every $y \in \Y$ we have an open neighborhood $\U$ around $\pi^{-1} (y)$ and $f \in \U$.  The conditions from the definition of Cartier divisors are satisfied. 
\qed

%The  morphism $U_n(\w) \rightarrow \prod_{i=1}^n \P_{\w_i, k}^1$ is called \textbf{weighted join construction}.
%For $\w=(1, \cdots, 1)$, we denote  by $\pi_n $ and $U_n$  instead of  $\pi_n (1, \cdots, 1)$ and
%$U_n(1, \cdots, 1)$.

%For each $i=1, \cdots, n$, we consider   $\P^1_{i, k} =\Proj (\Q[t_0, t_i])$ and  $\P^1_{w_i, k} =\Proj (\Q[x_0, x_i])$.
%The  finite surjective morphism  $g_i: \P^1_{i, k}  \to  \P^1_{w_i, k}$ given by
%$g_i([t_0: t_i])= [t_0^{q_0}: t_i^{q_i}]$ leads to  finite surjective morphism  
%\begin{equation}
%	\begin{split}
	%	g: \prod_{i=1}^n \P_{i, k}^1   &  \longrightarrow \prod_{i=1}^n \P_{\w_i, k}^1,  \\
	%		\x=[x_0: \cdots: x_n]  &\mapsto  \pi_{n, \w}:=( [x_0:x_1], \cdots,  [x_0:x_n]),
	%	\end{split}  
%\end{equation}

%As an application of  the above lemma we have the  following result.

\begin{prop}\label{p-blow}  
	% Silverman page 7 (top of page)  % reference [23] %Keeping the above notation,   
	Let $\X:= \prod_{i=1}^n \P_{\w_i, \Q}^1  $,  and consider   $\pi: \tX \rightarrow \X$,    the blow-up of $\X$ along  ${\bar 0}=([1:0], [1:0],\cdots,[1:0])$.
	Denote by  $\tY=\pi^{-1}({\bar 0})$   the exceptional divisor of this blow-up.    Then, for all $\nu\in M_\Q$ and any non-singular points 
	\[
	\x=[x_0: x_1 :\cdots: x_n] \in \P_{\w, \Q}^n \backslash \{ [1:0: \cdots: 0]\} 
	\]
	with $\bar x= \varphi_{n, \w}\mid_{\X} (\x) \in \X (\Q) \setminus \{{\bar 0}\}$, we have 
	\begin{equation}
		\il_{\tX, \tY  } ( \pi^{-1} (\bar x), \nu) 
		= \min \left\{ \left\lfloor  \frac {\nu^+ (x_0)}  {q_0}   \right\rfloor,  \dots ,  \left\lfloor  \frac {\nu^+ (x_n)}  {q_n}   \right\rfloor  \right\}=T_\nu(\bar x),
	\end{equation}
	and 
	% adding over all $\nu\in M_\Q$ gives that
	\begin{equation}
		\log \hwgcd (\x) = \sum_{\nu \in M_\Q}  \il_{\tX, \tY } ( \pi^{-1} (\bar x), \nu) = \lwh_{\tX, \tY} \left(\pi^{-1} (\bar x) , \nu \right).
	\end{equation}
\end{prop}

\proof    
By  Lem. \ref{blow-div},   $\tY$ is an effective divisor on $\tX$. Hence,    using the functoriality of local weighted heights, we have 
\begin{align*}
	\il_{\tX, \tY  } ( \pi^{-1} ({\bar x}), \nu)  & =  \il_{\X, {\bar 0} } ( { \bar x}, \nu) 
	%=\lambda_{X, {\bar 0}} ( { \bar x}, \nu) \\
	=   \il_{\P_{\w, \Q}^n, {[1: 0 : \cdots : 0]}}  \left( [x_0^{\frac{1}{q_0}}: x_1^{\frac{1}{q_1}} : \cdots : x_n^{\frac{1}{q_n}}] , \nu \right),  \\
	&= \min \left\lbrace  \nu^+ (x_0^{\frac{1}{q_1}}) , \cdots , \nu^+ (x_n^{\frac{1}{q_n}}) \right\rbrace \\
	& = \min \left\{ \left\lfloor  \frac{\nu^+ (x_0)}{q_0}     \right\rfloor,  
	\cdots ,  \left\lfloor  \frac {\nu^+ (x_n)}{q_n}    \right\rfloor  \right\}=T_\nu(\bar x).
\end{align*}
Adding these weighted local heights  together  we get the global formula.

\qed

The above result   leads  to the following definition.

\begin{defi}\label{glwgcd}
	Let $\X $ be a smooth weighted variety  defined over $\Q$, and $\Y $ a subvariety of $\X $ of codimension $r\geq 2$ and
	$\pi : \tX \to \X$, 
	the  blow-up of $\X$ along $\Y$.  For any $P \in (\X \setminus \Y)(\Q)$ we denote by 
	\[\tilde P := \pi^{-1} (P) \in \tX \quad \text{ and } \quad\tY=\pi^{-1}({\Y}).
	\] 
	The \textbf{generalized logarithmic weighted greatest common divisor of the point $P$ with respect to $\Y$ } is  defined to be 
	\begin{equation}
		\log \hwgcd (P; \Y) =  \lwh_{\tX, \tY} (\tilde P).
	\end{equation}
\end{defi}
%
%This is the \emph{weighted} analog of the generalized logarithmic   greatest common divisor defined in \cite{MR2162351}.

A point $\x=[x_0: \cdots: x_n] \in \P^n_{\w, \Q}$ is called 	 \emph{normalized} if it has integers coordinates and  $\wgcd (x_0, x_1, \cdots, x_n)=1$; see \cite{b-g-sh} for details.

\begin{lem} \label{peq0}
	Let $\w=(q_0, \cdots, q_n)$ be    well-formed weights,     $m=\lcm(q_0, \cdots, q_n)$, 
	$\y=[1: 0: \cdots: 0]$,  and     	  $\x=[x_0: \cdots: x_n] \in \P^n_{\w, \Q}$ a smooth and  normalized point. 
	%$ \x=[x_0: \cdots: x_n] \in \P^n_{\w, \Q} \backslash  \sing (\P^n_{\w, \Q}), $    such that   $x_i\in \Z$  and $\wgcd (x_0, x_1, \cdots, x_n)=1$.
	Then   
	\begin{equation}
		\label{estim}
		\log \hwgcd (\x;  \{\y \})=   \log \gcd (x_1, \dots, x_n) + O(1).
	\end{equation}
\end{lem}

\proof
Indeed,  letting $\w_i=(q_0, q_i)$ for each $i=1, \cdots, n$ and   
\[
\X=\prod_{i=0}^n \P_{\w_i, \Q}^1,
\]
then  considering the rational map      $ \pi_{n, \w}: \P_{\w , \Q }^n \rightarrow \X$, 
we have  ${\bar 0} =\varphi_{n, \w} (\y )$,  where  
\[
{\bar 0} =([1:0], \cdots, [1:0])  \in  \X.
\]
Let  ${\bar x}=\varphi_{n, \w} (\x )$  and  apply  Prop.   \ref{p-blow}, to the blow-up   $\pi: \tX \to \X$ along $\Y= \{\bar 0\}$.
Let $\tY=\pi^{-1}(\Y)$ be the exceptional  divisor of the blow-up.
Then
\[
\log \hwgcd (\x; \{\y \}) =\log \hwgcd (\bar{x}; \Y) =  \lwh_{\tX,\tY} \left( \pi^{-1} (\bar{x}) \right)
\]
By  definition of the global weighted height  and  properties   of local weighted height \cite[Thm. 1 (iv)]{Salami-Shaska1}, one can see that the last term is equal  to  the right-hand side of \eqref{estim}.
\qed
%

%^^^^^^^^^^^^^
One can extend the result of  Lem. \ref{peq0}, as stated in following proposition. 

%^^^^^^^^^
\begin{prop}\label{p-whp}
	Let $\w=(q_0, \cdots, q_n)$ be a   well-formed set of weights  and     $m=\lcm(q_0, \cdots, q_n)$. 
	Assume that $\Zc \subset \P_{\w, \Q}^n$  is a closed subvariety defined by the weighted homogeneous polynomials 
	$f_1, \cdots, f_r$  with integer coefficients such that 
	\[
	\Zc \cap  \sing(\P_{\w, \Q}^n )=\emptyset.
	\]
	Then 
	\begin{equation}
		\log \hwgcd (\x; \Zc)=    
		\log \gcd (f_1(\x ) , \dots, f_r(\x)) + O(1),
	\end{equation}
	for   $\x \in\P_{\w, \Q}^n \backslash  \left\lbrace \sing(\P_{\w, \Q}^n ) \cup \Supp(\Zc) \right\rbrace $   with $x_i\in \Z$ and $\wgcd(x_0, x_1, \cdots, x_n)=1.$
\end{prop}

\proof 
Let $\Y$ be given by
\begin{equation}\label{phi}
	\Y = \varphi_{n, \w} (\Zc) \subset \X=\prod_{i=1}^n \P_{\w_i, \Q}^1,
\end{equation}
where $\w_i=(q_0, q_i)$  for $i=1, \cdots, n$ and where $\varphi_{n, \w} $ is defined by \eqref{phi1}.

Consider the blow-up $\pi: \tX \to \X$ along $\Y$ and its exceptional divisor  $\tY=\pi^{-1} (\Y)$.
Let 
$\y=\varphi_{n, \w}(\x),$
for any 
$
\x \in\P_{\w, \Q}^n \backslash  \left\lbrace \sing(\P_{\w, \Q }^n ) \cup \Supp(\Zc) \right\rbrace.
$
Then,  we have 
\begin{align*}
	\log \hwgcd (\x ; \Zc)  & = \sum_{\nu_p \in M_\Q } 	\il_{\tX , \tY } ( \pi^{-1} (\y), \nu_p )  \\
	& = \sum_{\nu_p \in M_\Q }  \il_{\X , \Y }  ( \y, \nu_p )   \\
	&	 = \sum_{\nu_p \in M_\Q }  \il_{\P^n_{\w, \Q} , \Zc }  ( \x, \nu_p ).\\
	&  =\sum_{p \in M_\Q\backslash \{ \infty\} }   \il_{\P^n_{\w, \Q} , \Zc }  ( \x, \nu_p ) + 
	\il_{\P^n_{\w, \Q} , \Zc }  ( \x, \nu_\infty )   		
\end{align*}
By  definition (see \cite[Eq. (37)]{Salami-Shaska1}) of the local weighted height associated to $\Zc$, for   $\nu_p\in M_\Q$,  we have 
\begin{align*}
	\il_{\P^n_{\w, \Q} , \Zc }  ( \x, \nu_p )  & = \min_{1\leq j \leq r}^{} \left\lbrace   - \log \frac{|f_j(\x)|_{\nu_p}}{ \max_{i}^{} |x_i|_{\nu_p}^{\frac{m}{q_i}}}\right\rbrace \\
	& =  \log  \max_{1\leq j \leq r}^{}   |f_j(\x)|_{\nu_p} +  \frac{m}{q_i}\log  \max_{i}^{} |x_i|_{\nu_p}.
\end{align*}
Thus, 
\[
\sum_{\nu_p \in M_\Q  }   \il_{\P^n_{\w, \Q} , \Zc }  ( \x, \nu_p ) =  \log   \gcd (f_1(\x ) , \dots, f_t(\x)) +   O(1).
%  \max_{1\leq j \leq r}^{}   \log   |f_j(\x)| +     \frac{1}{q_i}\log  \max_{0\leq i \leq n}^{} |x_i|. 
\]
Therefore, putting all together gives the desired equality.
\[	
\log \hwgcd (\x ; \Zc) =   \log   \gcd (f_1(\x ) , \dots, f_t(\x)) + O(1).
\]
This completes the proof.
\qed

%^^^^^^^^^^^^^^^^^^^^^^*
\section{Vojta's conjecture and bounds on greatest common divisors}
Let $\X $ be a smooth weighted variety  defined over $\Q$, and $\Y $ a subvariety of $\X $ of codimension $r\geq 2$ and
$\pi : \tX \to \X$,  	the  blow-up of $\X$ along $\Y$.
It is a well-known fact that  the canonical bundle of the blow-up  $\pi: \tX \rightarrow \X$ along the subvariety
$\Y \subset \X$  is given by 
\[ 
K_{\tX} \sim \pi^\star  K_\X + (r-1) \tY,
\]
where  $\tY=\pi^{-1}(\Y)$.
If $\Ac$ is an ample divisor on $\X$, then there is an integer $N$ such that 
$
- \tY + N \pi^\star \Ac
$
is ample on $\tX$, see \cite[Thm. A.5.1]{Ha}.
Let  
\[
\tilde \Ac = - \frac 1 N \tY + \pi^\star \Ac
\]
be the ample cone of $\tX$.   We further assume that  $-K_\X$ is a normal crossing  and  
\[
\Supp (K_\X) \cap \Y= \emptyset.
\]
Let $S$ be a finite set of places of $\Q$ and  define 
\begin{equation}
	\lwh_{\X, D, S} (\cdot) := \sum_{\nu\in S} \il_{\X, D} (\cdot, \nu)   \quad \text{and } \quad 
	\lwh^\prime_{\X, D, S} (\cdot) :=  \sum_{\nu\not\in S} \il_{\X, D} (\cdot, \nu)   
\end{equation}
Then we have the following:

%^^^^^^^*

\begin{thm} \label{main0}
	Let $\X $ be a smooth projective weighted variety, $\Ac$ an ample divisor on $\X$, 
	$\Y \subset \X $ a smooth subvariety of codimension $r\geq 2$,
	and  $-K_\X$  a normal crossing divisor whose support does not intersect $\Y$, all defined over $\Q$.
	Assume     Vojta's conjecture (see Conj. \ref{voj1wsp}) for smooth weighted varieties. 
	Then for every finite set of places $S$ containing $M_\Q^\infty$ and every $0 < \varepsilon < r-1$ there is a proper closed subvariety 
	$ \Zc = \Zc (\varepsilon, \X, \Y, \Ac,   S) $ of $ \X,$
	and  constants  $C_\varepsilon = C_\varepsilon (\X, \Y, \Ac,   S)$ and  
	$\delta_{\varepsilon} = \delta_{\varepsilon} (\X, \Y, \Ac)$,  such that 
	\begin{equation}\label{meq1}
		\log \hwgcd (P; \Y) \leq \varepsilon \lwh_{\X, \Ac} (P) + \frac 1 {r-1+\delta_\varepsilon} \lwh^\prime_{\X, -K_\X, S} (P) + C_\varepsilon,
	\end{equation}
	for all $P\in (\X\setminus \Zc)  (\Q)$.
\end{thm}

\proof 
%The proof goes similarly to   \cite[Thm.~6]{MR2162351} with necessary adjustments.   
We apply    Conj. \ref{voj1wsp}, for the weighted  blow-up $\pi : \tX \rightarrow \X$ and the divisor $D =-\pi^\star K_\X$ to get 
\[
\lwh_{\tX, -\pi^\star K_\X, S} (\tP) + \lwh_{\tX, K_{\X}} (\tP) \leq \varepsilon \lwh_{\tX, \tilde \Ac} (\tP) + C_\varepsilon,
\]
for all $\tP \in (\tX  \setminus \tilde \Zc) (\Q)$.  
Substituting $K_{\tX}=\pi^\star K_\X + (r-1) \tY$ and $\tilde \Ac= - \frac 1 N \tY + \pi^\star \Ac$ we get 
\[
- \lwh_{\X, K_\X, S} (P) + \lwh_{\X, K_{\X}} (P) + (r-1) \lwh_{\tX, \tY} (\tP) \leq \varepsilon \lwh_{\X, \Ac}  (P) - \frac \varepsilon N \lwh_{\tX, \tY} (\tP)+ C_\varepsilon
\]
for all $P \in (\X   \setminus \pi(\tilde \Zc))(\Q)$.    Since
$ - \lwh_{\X, K_\X, S} (P) + \lwh_{\X, K_{\X}} (P)=\lwh^\prime_{\X, K_\X, S}
$
we have  
\[
\lwh^\prime_{\X, K_\X, S} (P)+ \left( r-1 + \frac \varepsilon N \right) \lwh_{\tX, \tY} (\tP)  \leq 
\varepsilon \lwh_{\X, \Ac} (P) + C_\varepsilon,
\]
for all $P \in (\X\setminus \Zc)  (\Q)$. Hence,
\[
\lwh_{\tX, \tY} (\tP)  \leq   \frac N {N(r-1) + \varepsilon} \left(  - \lwh^\prime_{\X, K_\X, S} (P)  + \varepsilon \lwh_{\tX, \Ac} (P) + C_\varepsilon \right)
\]
Since $\lwh_{\tX, \tY} (\tP) = \log \hwgcd (P; \Y)$, we have 
\[
\log \hwgcd (P; \Y) \leq   \frac N {N(r-1) + \varepsilon} \left(  - \lwh^\prime_{\X, K_\X, S} (P)  + \varepsilon \lwh_{\tX, \Ac} (P) + C_\varepsilon \right).
\]
Finally, setting $\delta=\varepsilon/N$ gives \eqref{meq1} and this completes the proof. 

\qed

%Using the above result we have  \cite[Thm.~2]{MR2162351}.
%^^^^^^^^^^^^^^^^^^^^^^^ 
%\subsection{Vojta's conjecture and gcds}
%
Let  $\w=(q_0, \cdots, q_n)$ and assume that  $X=\P_{\w, \Q}^n$ is well-formed weighted projective space. For    $i=0, \ldots , n,$ 
let 
$$H_i=\{ \x\in X \mid  x_i=0\}, $$
$A_0= H_0$,  and   $ K_X= - \sum_{i=0}^n H_i$. 
the map 
\[
\varphi_{n, \w}: \P^n_{\w, \Q} \to  \X :=\prod_{i=0}^n \P_{\w_i, \Q}^1,
\]
where $\w_i=(q_0, q_i)$  for $i=1, \cdots, n$ and  denote by
\[
\cH_i =\varphi_{n, \w}(H_i),   \quad       \Ac_0=\varphi_{n, \w}(H_0),  \quad    K_\X= - \sum_{i=0}^n \varphi_{n, \w}(H_i).
\]
Notice that the canonical divisor $K_\X$ is a normal crossing   on  $\X$ satisfying  
\[
\Y \cap \Supp (-K_\X)= \emptyset, \quad \text{ where }\Y=\varphi_{n, \w}(\Zc).
\]
%}

Let  $S$ be a   finite set of  primes. 
The ``prime-to-$S$'' part $|x|'_S$ of any nonzero integer $x$  is defined to be the largest divisor of $x$ that is not divisible by any of the primes in $S$, in other words, 
\[
|x|'_S=|x| \cdot \prod_{p \in S}^{} |x|_p.
\]
In particular, $x$ is an $S$-unit if and only if  $|x|'_S=1.$

\begin{thm}\label{main1}
Keeping the above notations in mine, we assume   that
$\Zc \subset \P_{\w, \Q}^n$  is a closed subvariety defined by
\[
f_1, \cdots, f_t \in \Z_\w [x_0, \ldots , x_n],
\] 
of codimension    $r:=n- \dim(\Z) \geq 2$ in $\P_{\w, \Q}^n$
such that $\Zc \cap  \sing(\P_{\w, \Q}^n )=\emptyset $, and has
transversal intersection with the
union  $\displaystyle \cup_{i=0}^n \cH_i$.
%} 
Let  $S$ be a   finite set of  primes and    $\varepsilon >0$.   
If     Conj. \ref{voj1wsp}, holds for the weighted blow-up $\pi: \tX \rightarrow \X$ along $\Y$,
then   there exists a nonzero weighted  polynomial 
$g \in \Z_\w [ x_0, \ldots , x_n]$ 
and a constant $\delta=\delta_{\varepsilon, \Zc} >0$,     such that every $\tilde \a= (\a_0, \cdots, \a_n) \in \Z^{n+1}$ with    $\wgcd(\a_0, \cdots, \a_n)=1$  satisfies either $g(\tilde \a)=0$ or
\begin{equation}\label{mm10}
\gcd(f_1({\tilde \a}), \cdots, f_t({\tilde \a}))\leq  \max \, \left\lbrace  |\a_0|^{\frac{1}{q_0}}, \cdots,  |\a_n|^{\frac{1}{q_n}}\}\right\rbrace^\varepsilon  \cdot \left(  |\a_0  \cdots \a_n|'_S\right)^{\frac{1}{ m(r-1+ \delta  )}},
\end{equation}
where $|\cdot|'_S$ is the ``prime-to-$S$'' part of its origin.
%, see \cite[Def. 1] {MR2162351}
\end{thm}

\proof        
By definition of the global weighted height  for  
\[
\x=[\a_0: \cdots: \a_n] \in \P^n_{\w, \Q}    \backslash \{ \sing(\P^n_{\w, \Q}) \cup \Supp (\Zc) \}(\Q)
\]
with $\wgcd(\a_0, \cdots, \a_n)=1$,   we have
\begin{equation}
\label{peq1}
\lwh_{X, A} (\x)= \log \, \max \{ |\a_0|^{\frac{1}{q_0}}, \cdots, |\a_n|^{\frac{1}{q_n}} \} + O(1)
\end{equation}
By  Prop.  \ref{p-whp},   we have 
\begin{equation}
\label{peq2}
\log \hwgcd \left(  \x ; \Zc \right) =    \log \gcd (f_1(\x), \cdots, f_t(\x)) + O(1).
\end{equation}
Let  $\y=\varphi_{n, \w}(\x)$   for  $\x=[\a_0: \cdots: \a_n] \in \P^n_{\w, \Q}$. Then,   
by definition of $S$-part of the weighted heights and functoriality of the weighted heights, we have
\[
\lwh^\prime_{\X, \cH_i, S}(\y)=\lwh^\prime_{X, H_i, S}(\x)= \sum_{\nu \in S}^{} \nu^+(\a_i)= \frac{1}{q_i} \log |\a_i|'_S,
\]
which implies that
\begin{equation}
\label{peq3}
\lwh^\prime_{\X,- K_\X, S} (  \y)=
\lwh^\prime_{X,- K_X, S} (\x)=\sum_{i=0}^n  \lwh^\prime_{X, H_i, S}= \frac{1}{m} \log |\a_0 \a_1 \cdots \a_n|'_S.
\end{equation}
By substituting \eqref{peq1},  \eqref{peq2}, \eqref{peq3}, in  \eqref{meq1},  one obtains 
\begin{align*}
\log \gcd (f_1({\tilde \a}), \cdots, f_t({\tilde \a})) & 
\leq    \log \hwgcd \left(  \x ; \Zc \right)=  \log \hwgcd \left(  \y ; \Y \right) \\
&  \leq \varepsilon    \lwh_{\X, \Ac} (\y) + \frac 1 {r-1+\delta} \lwh^\prime_{\X, -K_\X, S} (\y) + C_\varepsilon, \\
& \leq  \varepsilon \cdot \log \max  \{ |\a_0|^{\frac{1}{q_0}}, \cdots, |\a_n|^{\frac{1}{q_n}} \} \\
& \quad \ \ \  \cdot {\frac{1}{ m (r-1+ \delta)}} \log \left( |\a_0 \a_1 \cdots \a_n |'_S\right) + C_\varepsilon,
\end{align*}
where $\delta=\delta_{\varepsilon, \Zc}$.
%Replacing $\varepsilon$ with $\varepsilon/m$, then  multiplying  the both sides by $m$ and  
Then, by   exponentiation, we obtain   the desired inequality \eqref{mm1}. 
\qed

%\section{Acknowledgment} 

%We would like to thank the anonymous referees for their valuable comments and suggestions, which helped us significantly improve the manuscript.

\bibliographystyle{alpha} 
\bibliography{refs}% common bib file

\end{document}